\documentclass[letterpaper,11pt]{article}

\usepackage[T1]{fontenc}
\usepackage[utf8]{inputenc}
\usepackage{lmodern}

\usepackage[letterpaper, top=1in, bottom=1.1in]{geometry}
\usepackage{siunitx}
\usepackage{rotating}
\usepackage{hyperref}
\usepackage{authblk}
\usepackage[all]{nowidow}
\usepackage[format=plain, 
            labelfont={bf,it, footnotesize},
            textfont={it, footnotesize}]{caption}

\usepackage{amsthm, amsmath, amsfonts, amssymb}
\usepackage{mathtools}
\usepackage{bm}

\newtheorem{theorem}{Theorem}[section]

\newtheorem{corollary}{Corollary}
\newtheorem{conjecture}{Conjecture}

\DeclareMathOperator{\supp}{supp}
\DeclareMathOperator*{\argmax}{arg\,max}
\DeclareMathOperator*{\esssup}{ess\,sup}

\newcommand{\Ebb}{\mathbb{E}}
\newcommand{\Mbb}{\mathbb{M}}
\newcommand{\Nbb}{\mathbb{N}}
\newcommand{\Pbb}{\mathbb{P}}
\newcommand{\Qbb}{\mathbb{Q}}
\newcommand{\Rbb}{\mathbb{R}}

\newcommand{\Fcal}{\mathcal{F}}
\newcommand{\Hcal}{\mathcal{H}}
\newcommand{\Mcal}{\mathcal{M}}
\newcommand{\Pcal}{\mathcal{P}}
\newcommand{\Qcal}{\mathcal{Q}}
\newcommand{\Xcal}{\mathcal{X}}
\newcommand{\Zcal}{\mathcal{Z}}

\title{Fading Boundaries: On a Nonparametric Variant\\ of the Kiefer--Weiss Problem}
\author{Michael Fau{\ss} and H.~Vincent Poor}
\affil{\normalsize Department of Electrical Engineering, Princeton University\\
       \{mfauss, poor\}@princeton.edu}
\date{October 2020}

\begin{document}

\maketitle

\begin{abstract}
  A nonparametric variant of the Kiefer--Weiss problem is proposed and investigated. In analogy to the classical Kiefer--Weiss problem, the objective is to minimize the maximum expected sample size of a sequential test. However, instead of taking the maximum over a parametric family of distributions, it is taken over all distributions defined on the given sample space. Two optimality conditions are stated, one necessary and one sufficient. The latter is based on existing results on a more general minimax problem in sequential detection. These results are specialized and made explicit in this paper. It is shown that the nonparametric Kiefer--Weiss test is distinctly different from its parametric counterpart and admits non-standard, arguably counterintuitive properties. In particular, it can be nontruncated and critically depends on its stopping rules being randomized. These properties are illustrated numerically using the example of coin flipping, that is, testing the success probability of a Bernoulli random variable.
\end{abstract}

\section*{Prologue}  
\label{sec:prologue} 

Beenish the blind magician had two coins that she used for some of the tricks in her magic show. One of the coins was more likely to show heads, the other was more likely to show tails; otherwise the coins were indistinguishable. After each show, Beenish liked to ask her assistant, Sabri, to help her tell the coins apart. To this end, Beenish would take one of the coins, flip it, and ask Sabri to tell her whether it shows heads of tails. This procedure would then repeat until Beenish was sufficiently certain which coin she was flipping.

Over time, Sabri came to enjoy this little post-show ritual, which gave him a welcome break from feeding rabbits and polishing wands. At some point, he started to deliberately prolong the procedure by occasionally misreporting the outcome of a flip. This made him wonder just how long he could keep Beenish flipping, without her reaching a decision or getting tired of his shenanigans. 

Beenish had already figured out what Sabri was up to. However, instead of confronting him, she was curious to see if she could come up with a way of running their ritual that would prevent her from spending whole afternoons flipping coins, yet would reliably tell the coins apart on days when Sabri was reporting the true outcomes---which he kept doing most of the time.

In a first attempt at prolonging the ritual, Sabri tried reporting heads and tails alternatingly, which, he thought, should make it impossible for Beenish to come to a decision. For the first couple of flips, this strategy seemed to work well. However, as Beenish kept flipping, she grew increasingly skeptical of the repeating pattern and, much to Sabri's disappointment, finally stopped even sooner than usual. Sabri's second attempt was to first report several heads in a row, only to report a similar number of tails afterwards. However, Beenish soon started to recognize the same alternating pattern in these sequences. 

Not satisfied with his attempts thus far, Sabri came to the conclusion that he had to give up on using predictable patterns and instead had to emulate the randomness of the coins. He started bringing a regular coin with him that he secretly flipped whenever Beenish flipped her coin. This attempt turned out to be much more successful and resulted in a record number of flips. However, Beenish, knowing that she could no longer rely on simply looking out for patterns, soon changed her strategy. She had usually continued flipping until she first observed a certain difference in the numbers of heads and tails. She now chose this difference larger than she used to, but gradually reduced it as the ritual went on. In this way, she could still trust her decisions on days when the ritual ended early; yet, on days when it was dragging on, she could trade-off her confidence and the number of flips in a reasonable manner.

Despite Beenish's new strategy, Sabri was quite happy with what he had achieved. Even on days when he was telling the truth, the ritual now took slightly longer than before. However, on certain occasions he got frustrated with his fair coin giving him a ``wrong'' outcome, one that he knew would make Beenish stop the ritual. He did not want to fall into the trap of introducing predictable patterns again, yet he sought to have just a little more control over his coin. The next day, he passed by The Magic Depot and got a whole box full of coins, all with different probabilities of showing heads or tails. Now all that was left to do, was for him to figure out when to use which coin\ldots

Beenish soon started to notice that their coin flipping rituals had been taking longer recently. In particular, she noticed that on more occasions she would end up stopping the ritual simply because she had run out of patience---this had happened before, but very infrequently. After a particularly long Wednesday afternoon session, she decided to put an end to all of this. From now on, she would decide on the number of flips beforehand, so Sabri could fiddle with his coins as much as he liked. However, after some particularly boring and lifeless rituals of this kind, she decided that this could not be the end of it. She still had a couple of tricks up her sleeve and a strong conviction that there was a way to keep the ritual fun for Sabri, while at the same time giving her thumb a well-deserved rest.

\section{Introduction}   
\label{sec:introduction} 

The type of problem illustrated in the Prologue is known as the \emph{Kiefer--Weiss problem} in sequential analysis, named after Jack Kiefer and Lionel Weiss who first proposed and studied it in the 1950s \cite{Kiefer1957, Weiss1962}. 

The observation that led to the formulation of the Kiefer--Weiss problem is that sequential hypothesis tests, while being highly efficient under ideal conditions, show undesirable behavior under model mismatch. In particular, their expected sample sizes can increase significantly if the true distribution is equally similar to both hypothesized distributions---a classic example being that of a sequential test for a mean parameter being positive or negative, when the true mean is zero. Since a reduced expected sample size is often the main reason for using sequential tests in the first place, this effect can be highly problematic in delay-critical applications; a more detailed discussion of this aspect is deferred to the next section. 

Kiefer and Weiss proposed to design a sequential test such that, in addition to it meeting the targeted error probabilities under both hypotheses, it minimizes the \emph{maximum} expected sample size over a parametric family of distributions. This problem, in various varieties, has received considerable attention in the literature \cite{Dvoretzky1953, Kiefer1957, Anderson1960, Weiss1962, Robbins1972, Lai1973, Lorden1976, Eisenberg1982, Huffman1983, Dragalin1988, Lai1988, Pavlov1991, Augustin2001, Zhitlukhin2013}. 

In this paper, a variation of the classical Kiefer--Weiss problem is studied. Again, the objective is to minimize the maximum expected sample size, however, the maximum is taken over \emph{all} possible distributions on the given sample space, instead of over a parametric family. Hence, the problem investigated in this paper is referred to as the \emph{nonparametric} Kiefer--Weiss problem.

In principle, the nonparametric Kiefer--Weiss problem is a special case of the minimax sequential hypothesis testing problem studied in \cite{Fauss2020_aos}. As a consequence, most of the findings presented here are based, either directly or indirectly, on results in \cite{Fauss2020_aos}. Nevertheless, we consider the nonparametric Kiefer--Weiss problem to be an interesting special case in its own right. First, it is simple enough to make some of the implicit optimality conditions given in \cite{Fauss2020_aos} explicit, which allows for a more tangible characterization of the optimal test. Yet, second, it is complex enough to lead to an interesting, non-standard optimal test whose properties are notably different from typical sequential procedures. Third, the Kiefer--Weiss problem is not only interesting from a theoretical point of view, but also relevant from a practical and historical perspective. This aspect is discussed in more detail in Section~\ref{sec:model_mismatch}.

A formal statement of the nonparametric Kiefer--Weiss problem is given in Section~\ref{sec:nonpa_kiwei}, where also two optimality conditions for the corresponding optimal test are stated, one being necessary and one being sufficient. The properties of a test satisfying the sufficient optimality conditions are then discussed in Section~\ref{sec:nonpa_kiwei_properties}. An illustrative example of a nonparametric Kiefer--Weiss test, based on the coin flipping scenario in the Prologue, is given in Section~\ref{sec:example}. Section~\ref{sec:conclusion} concludes the paper and provides a brief outlook on open questions and possible future work, including a conjecture on a connection to fixed-sample-size tests.

\subsection{Notation and Assumptions} %
Throughout the paper, $\Nbb$ ($\Nbb_0$) denotes the the set of positive (nonnegative) integers and $\Rbb_+^K$ denotes the space of $K$-dimensional nonnegative real-valued vectors. No distinction is made between row and column vectors. The element-wise product of two vectors $\bm{x}$ and $\bm{y}$ is denoted by $\bm{xy}$. All comparisons are defined element-wise. The notation $\partial_{z_k} f(\bm{z})$ is used for the superdifferential \cite[\S 23]{Rockafellar1970} of a concave function $f \colon \Rbb_+^K \to \Rbb_+$ with respect to (w.r.t.) $z_k$ evaluated at $\bm{z}$. Note that $\partial_{z_k} f(\bm{z})$ is a left-closed interval on the real line. For concave functions, the minimum of this interval coincides with the right-derivative and is denoted by $\min \partial_{z_k} f(\bm{z}) = \partial_{z_k}^+f(\bm{z})$. Note that at $z_k = 0$ the superdifferential is unique and coincides with the right-derivative, that is, $\partial_{z_k} f(\bm{z}) |_{z_k = 0} = \{ \partial_{z_k}^+f(\bm{z}) |_{z_k = 0}\}$. If a function $h$ exists such that $h(\bm{z}) \in \partial_{z_k} f(\bm{z}) \; \forall \bm{z} \in \Rbb_+^K$, then $h$ is called a partial derivative of $f$ w.r.t.\ $z_k$. The set of all partial derivatives is denoted by $\partial_{z_k} f$.

In what follows, $\bm{X} = X_1, X_2, \ldots$ denotes a sequence of random variables taking values in a measurable space $(\Xcal, \Fcal)$. The subsequence $(X_1, \ldots, X_n)$, $n \geq 1$, is denoted by $\bm{X}_n$. Realizations of $\bm{X}$ and $\bm{X}_n$ are denoted by $\bm{x}$ and $\bm{x}_n$, respectively. The joint distribution of the sequence $\bm{X}$ is denoted by $\Pbb$. The shorthand $\Pbb$-a.e. is used for properties that hold $\Pbb$-almost everywhere. The notation $\Pbb = P^\Nbb$, where $P$ is a distribution on $(\Xcal, \Fcal)$, is used to indicate that all elements of $\bm{X}$ are independent and identically distributed (i.i.d.). This notation has to be read in the sense that the joint distribution of $\bm{X}_n$ under $\Pbb$ is given by $P^n$ for all $n \in \Nbb$. The set of all $P$ that admit a density $p$ w.r.t.\ a suitable background measure $\mu$ is denoted by $\Mcal_\mu$. Analogously, the set of distributions admitting a density w.r.t.\ $\mu^\Nbb$ is denoted by $\Mbb_\mu$.

The two hypotheses about the true distribution $\Pbb$ are written as 
\begin{equation}
  \begin{aligned}
    \Hcal_1 \colon \Pbb &= \Pbb_1, \\ 
    \Hcal_2 \colon \Pbb &= \Pbb_2.
  \end{aligned}
  \label{eq:hypotheses}
\end{equation}
In the literature, denoting the hypotheses by $\Hcal_0$ and $\Hcal_1$ is more common and, admittedly, makes it easier to identify a null hypotheses. However, as will become clear in the course of the paper, the notation in \eqref{eq:hypotheses} is more appropriate in the context of the Kiefer--Weiss problem, where the distinction between null and alternative hypothesis is of little importance. Moreover, since many of the results presented here are based on results in \cite{Fauss2020_aos}, it is helpful to use the same notation.

Both distributions $\Pbb_1$ and $\Pbb_2$ are assumed to be product distributions of the form $\Pbb_1 = P_1^\Nbb$ and $\Pbb_2 = P_2^\Nbb$ for some distributions $P_1, P_2 \in \Mcal_\mu$, that is, $\bm{X}$ is assumed to be an i.i.d.~process under both $\Hcal_1$ and $\Hcal_2$.\footnote{
  An extension to processes with Markovian representations, which is the scenario studied in \cite{Fauss2020_aos}, should be possible, but will not be attempted here.
}
Naturally, it is assumed that $P_1 \neq P_2$, more precisely that $D_\text{KL}(P_1 \Vert P_2) > 0$, where $D_\text{KL}$ denotes the Kullback--Leibler divergence.

In general, a sequential test for two hypotheses is specified via two sequences of randomized decision rules, $\psi = \bigl( \psi_n \bigr)_{n \geq 0}$ and $\delta = \bigl( \delta_n \bigr)_{n \geq 0}$. Each $\psi_n \colon \Xcal^n \to [0,1]$ denotes the probability of stopping the test after the $n$th sample has been observed, and each $\delta_n \colon \Xcal^n \to [0,1]$ denotes the probability of deciding for $\Hcal_1$, given that the test has stopped. The randomization is assumed to be performed by independently drawing from a Bernoulli distribution with success probability $\psi_n$ or $\delta_n$, respectively. The set of randomized decision rules defined on $\Xcal^n$ is denoted by $\Delta_n$. For the sake of a more concise notation, let $\pi = \bigl( \pi_n \bigr)_{n \geq 0}$, with $\pi_n = (\psi_n, \delta_n) \in \Delta_n^2$, denote a sequence of tuples of stopping and decision rules. In what follows, $\pi$ is referred to as a \emph{testing policy} and the set of all feasible policies is denoted by $\Pi \coloneqq \bigtimes_{n \geq 0} \Delta_n^2$.

The stopping time of a test with policy $\pi$, given the realization $\bm{X} = \bm{x}$, is denoted by $\tau_\pi(\bm{x})$. In order to define $\tau_\pi$ formally, let $\bigl( B_n \bigr)_{n \geq 0}$ be a sequence of independent Bernoulli random variables with success probabilities $\psi_n(\bm{x})$. The stopping time is then given by
\begin{equation}
  \tau_\pi(\bm{x}) = \min\{\, n \geq 0 : b_n = 1 \,\},
  \label{eq:stopping_time}
\end{equation}
where $b_n$ denotes the realization of $B_n$. The expected sample size of a test using policy $\pi$ under distribution $\Pbb$ is defined as
\begin{equation}
  \gamma(\pi,\Pbb) \coloneqq \Ebb_{\pi,\Pbb}\bigl[\, \tau_\pi(\bm{X}) \,\bigr], 
  \label{eq:runlength}
\end{equation}
where $\Ebb_{\pi,\Pbb}$ denotes the joint expectation w.r.t.\ the distribution $\Pbb$ of $\bm{X}$ and the randomization of the policy $\pi$. Analogously, the error probabilities of a test using policy $\pi$ are given by
\begin{align}
  \alpha_1(\pi,\Pbb_1) &\coloneqq \mathbb{E}_{\pi,\Pbb_1}\bigl[\, 1-\delta_{\tau_\pi}(\bm{X}) \,\bigr], \\ 
  \alpha_2(\pi,\Pbb_2) &\coloneqq \mathbb{E}_{\pi,\Pbb_2}\bigl[\, \delta_{\tau_\pi}(\bm{X}) \,\bigr].
  \label{eq:errors}
\end{align}
The expected sample size and the error probabilities are the performance measures used in the formulation of the Kiefer--Weiss problem, with the expected sample size being of particular interest. An more in-depth introduction to the Kiefer--Weiss problem and its rationale are given in the next section.

\section{Model Mismatch and Sequential Tests} 
\label{sec:model_mismatch}                    

In order to understand the motivation for the nonparametric Kiefer--Weiss problem, it is instructive to trace its origins, starting with Wald's famous sequential probability ratio test (SPRT) and its behavior under model mismatch. Readers who are already familiar with the Kiefer--Weiss problem or are not interested in this discussion can directly skip to Section~\ref{sec:nonpa_kiwei}.

\subsection{The Sequential Probability Ratio Test} %
Let $\bigl( P_\theta\bigr)_{\theta \in \Theta}$ be a parametric family of distributions on $(\Xcal, \Fcal)$ and let $\Pbb_\theta = P_\theta^\Nbb$. Wald's sequential detection problem is given by \cite{Wald1947, Poor2008}
\begin{equation}
  \begin{aligned}
    &\min_{\pi \in \Pi} \; \gamma(\pi, \mathbb{P}_{\theta_0}) & &\text{s.~t.} & \alpha_1(\pi, \mathbb{P}_{\theta_1}) &\leq \overline{\alpha}_1, \\
   & & & & \alpha_2(\pi, \mathbb{P}_{\theta_2}) &\leq \overline{\alpha}_2,
  \end{aligned}
  \label{eq:sprt}
\end{equation}
where $\theta_0 \in \{ \theta_1, \theta_2 \}$ and $\overline{\alpha}_1, \overline{\alpha}_2 \in (0, 1)$. That is, the expected sample size under $\Hcal_1$ or $\Hcal_2$ is minimized, subject to constraints on both error probabilities. Alternatively, the problem in \eqref{eq:sprt} can be formulated as an unconstrained optimization, where the error probabilities are included in the cost function:
\begin{equation}
  \min_{\pi \in \Pi} \; \gamma(\pi, \mathbb{P}_{\theta_0}) + \lambda_1 \alpha_1(\pi, \mathbb{P}_{\theta_1}) + \lambda_2 \alpha_2(\pi, \mathbb{P}_{\theta_2}).
  \label{eq:sprt_unconstrained}
\end{equation}
In \eqref{eq:sprt_unconstrained}, $\lambda = (\lambda_1, \lambda_2)$ are two positive cost coefficients that need to be chosen such that the test admits the desired error probabilities. In what follows, the focus will be on the unconstrained problem formulation in \eqref{eq:sprt_unconstrained} since, in the context of this paper, is easier to work with while proving the same conceptual insights.

Wald and Wolfowitz \cite{Wald1948} showed that the problem in \eqref{eq:sprt_unconstrained} is solved by the SPRT, which is a likelihood ratio test with two constant thresholds; compare Figure\ref{fig:thresholds}. While the SPRT is striking in its simplicity and has been immensely successful and influential \cite{Lai2001, Todd2007}, it does come with its own drawbacks. From the beginning, it was pointed out by Wald \cite{Wald1947} and others \cite{Anderson1960, Bechhofer1960} that the increase in sample efficiency depends on how much the true distribution deviates from the assumed ones. If the model mismatch becomes too large, a sequential test can in fact require \emph{more} samples on average than its fixed-sample-size counterpart.

\begin{figure}
  \centering
  \includegraphics{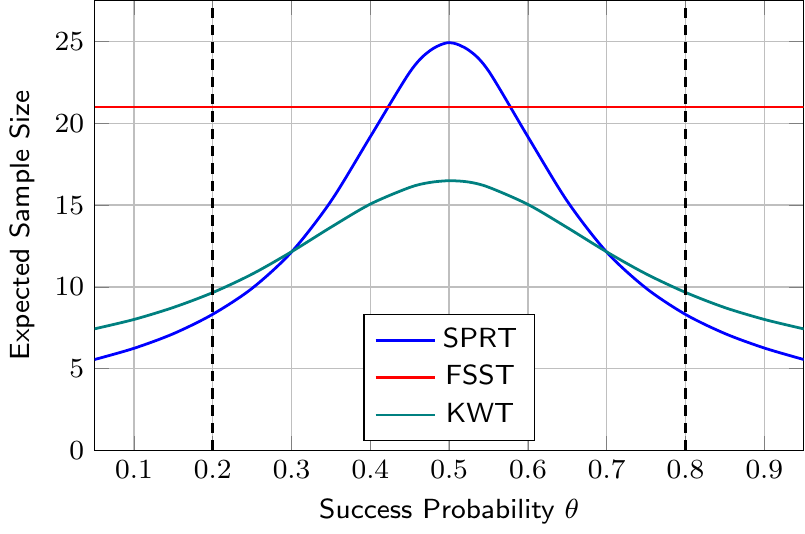}%
  \caption{Expected sample size of a sequential probability ratio test (SPRT), a fixed-sample-size test (FSST) and a Kiefer--Weiss test (KWT) for the hypotheses in \eqref{eq:hyp_illustration}. All tests were designed such they admit approximately identical error probabilities of $\alpha_1 = \alpha_2 \approx 10^{-4}$.}
  \label{fig:sample_size}
\end{figure}

This effect can be illustrated using the example from the Prologue. Let $\theta \in (0,1)$ denote the success-probability parameter of a Bernoulli distribution with probability mass function (PMF)
\begin{equation}
  p_\theta(x) = \theta^x (1-\theta)^x, \quad x \in \{0, 1\},
  \label{eq:nominal_densities}
\end{equation}
and consider a test for the two simple hypotheses
\begin{equation}
  \begin{aligned}
    \Hcal_1 &\colon \theta = \theta_0 = 0.8, \\
    \Hcal_2 &\colon \theta = \theta_1 = 0.2.
  \end{aligned}
  \label{eq:hyp_illustration}
\end{equation}
The blue line in Figure~\ref{fig:sample_size} shows the expected sample size of the corresponding SPRT as a function of the true success probability $\theta$. The thresholds were chosen such that $\alpha_1 = \alpha_2 \approx 10^{-4}$. For comparison, the number of samples required by a fixed-sample-size test (FSST) with comparable error probabilities ($\alpha_1 = \alpha_2 \approx 9.2 \times 10^{-5}$) is plotted in red. At the hypothesized parameter values, the higher efficiency of the SPRT is clearly visible, its expected sample size being only about one third of that of the FSST. However, it can also be seen that around $\theta = 0.5$ the expected sample size the SPRT exceeds that of the FSST by about four samples.

But why is this effect problematic in the first place? It is clear that neither the SPRT nor the FSST provide meaningful results if the true parameter value cannot be associated with either hypothesis. In fact, for $\theta = 0.5$ in the above example, both tests are equally likely to accept or reject $\Hcal_1$ so that the outcome of the test provides no information. However, there are two critical issues that only arise with the SPRT. First, while on average it might only require a few samples more%
\footnote{
  The Bernoulli distribution is far from being a worst-case scenario. For other distributions the increase in expected sample size can be much more severe; compare, for example, Figs.~3.4 and 5.1 in \cite{Tartakovsky2014}.
}
than the FSST, the right tail of the stopping time distribution can be shown to decay (only) exponentially, so that significantly larger stopping times occur with non-negligible probability. In this example, the sample size of the SPRT at $\theta = 0.5$ exceeds \num{65} in approximately \SI{5}{\percent} of cases. Such exceptionally large delays are often unacceptable, in particular if the sequential test is just one element of a larger processing chain. Also, the sensitivity of the SPRT to model mismatch can lead to situations where a system works well under regular conditions, but becomes unresponsive under irregular conditions. This can be fatal when the irregularity is caused by some underlying system- or sensor-failure which in turn calls for a quick intervention. 

These considerations point towards a second issue with sequential tests under model mismatch, namely, that they do not provide a reliable indicator for it. For a FSST, the value of the test statistic provides additional information about how trustworthy the corresponding decision is, typically in form of a p-value. A test statistic barely exceeding the threshold is much less indicative than a test statistic exceeding the threshold by a large margin. This confidence information gets lost when using the SPRT, whose test statistic is \emph{always} close to one of the thresholds. Hence, the only indicator for a mismatch is the sample size itself, which is much less informative than a p-value. In combination, these issues can prohibit the use of standard SPRTs in many applications. 

A straightforward remedy, which was already proposed by Wald \cite{Wald1947}, is to stop the test after a given maximum number of samples, irrespective of whether a threshold has been crossed. Tests of this type are known as \emph{truncated} sequential tests \cite{Genizi1965, Tantaratana1977, Tantaratana1982, Tartakovsky2014}. In practice, truncated sequential tests are usually the method of choice since they combine the predictable worst case behavior of FSSTs with the better performance under nominal conditions of the SPRT. Moreover, unless the test is truncated very aggressively, the thresholds of the standard SPRT can be used virtually unaltered, in particular if they were chosen based on Wald's approximation \cite{Wald1947}, which provides a conservative upper bound in the first place. More elaborate ways of compensating for the truncation can be found in the literature \cite{Bussgang1967, Madsen1974, Hu2011}. 

While the truncated SPRT is an excellent choice in practice, it is arguably less satisfying from a theoretical point of view. In particular, it lacks the property of both the SPRT and the FSST to be the solution of a well-defined optimization problem. In other words, it is not clear which \emph{design objective} underpins the truncated SPRT, apart from the notion of ``making the SPRT more robust against model mismatch.'' In order to close this gap, Kiefer and Weiss proposed a more principled approach to fixing the shortcomings of the SPRT. 

\subsection{The Kiefer--Weiss Test} %
The Kiefer--Weiss problem is formulated in close analogy to the sequential testing problem in \eqref{eq:sprt_unconstrained}, with the difference that it aims at minimizing the expected sample size under the least favorable parameter value \cite{Weiss1962, Tartakovsky2014}:
\begin{equation}
    \min_{\pi \in \Pi} \; \sup_{\theta_0 \in \Theta} \; \gamma(\pi, \mathbb{P}_{\theta_0}) + \lambda_1 \alpha_1(\pi, \mathbb{P}_{\theta_1}) +  \lambda_2 \alpha_2(\pi, \mathbb{P}_{\theta_2}).
  \label{eq:kiwei}
\end{equation}
While this is clearly not the only way of formalizing the notion of a sequential test with well-behaved sample sizes under model mismatch, it sits well in the general framework of minimax robustness \cite{VerduPoor1984, Minimax1995} and has attracted considerable attention in the literature. 

The general solution of the Kiefer--Weiss problem turned out to be rather elusive, which caused the focus to shift to special cases and approximations. The most well-studied cases are the Kiefer--Weiss test (KWT) for the mean parameter of a normal distribution \cite{Anderson1960, Weiss1962, Lai1973, Lorden1976} and the closely related problem of testing the drift parameter of a Wiener process \cite{Dvoretzky1953, Zhitlukhin2013}. Asymptotic and approximate results exist for more general families of distributions \cite{Robbins1972, Eisenberg1982, Huffman1983, Dragalin1988, Pavlov1991, Augustin2001}. An up-to-date treatment of the Kiefer--Weiss problem and some of its variations can be found in \cite[Section~5.3]{Tartakovsky2014}.

\begin{figure}
  \centering
  \includegraphics{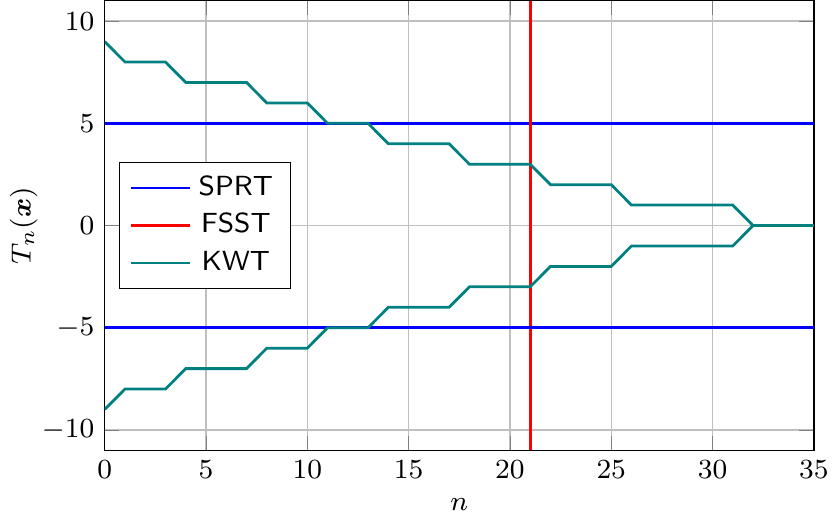}
  \caption{Thresholds of a sequential probability ratio test (SPRT), a fixed-sample-size test (FSST) and a Kiefer--Weiss test (KWT) for the hypotheses in \eqref{eq:hyp_illustration}. The test statistic on the ordinate is the difference of successes and failures. All tests were designed such they admit approximately identical error probabilities of $\alpha_1 = \alpha_2 \approx 10^{-4}$.}
  \label{fig:thresholds}
\end{figure}

One of the special cases of the Kiefer--Weiss problem that can be solved exactly is a test for the success-probability parameter of a binomial distribution \cite{Lorden1976}, which includes the coin-flipping example introduced above. For $\alpha_1 = \alpha_2 \approx 10^{-4}$, the thresholds of a KWT for the hypotheses in \eqref{eq:hyp_illustration} are given by the teal line in Figure~\ref{fig:thresholds}. For comparison, the thresholds of the regular SPRT and the FSST are depicted as well. Note that the thresholds are given in terms of the difference of the number of successes and failures instead of the likelihood ratio, that is, the test statistic is given by
\begin{equation}
  T_n(\bm{x}) = 
  2 \sum_{i=1}^n x_i - n.
\end{equation}

In general, the KWT is of the same form as the SPRT in the sense that it compares the likelihood ratio of the current sample to a pair of thresholds. However, the thresholds of the KWT depend on the sample number. Initially, they are larger than those of the SPRT, but then decrease as the sample size increases. Finally, both thresholds intersect, thus truncating the test. This ``soft'' truncation is different from the ``hard'' truncation of the trucated SPRT since it avoids situations in which a decision has to be made without the test statistic having crossed a threshold. Moreover, the successive reduction in confidence required by the KWT leads to an interesting coupling between the test statistic and the sample size. While the SPRT provides information about a possible model mismatch only in terms of the sample size, and the FSST only in terms of a confidence value, the KWT provides a combination of both, with larger sample sizes implying lower confidence values and vice versa.

In practice, however, the KWT is rather unpopular, occasional exceptions aside \cite{Mulekar1992}. The reason for this is that the ``optimal compromise'' between nominal and worst-case sample size that it offers comes at the cost of an increase in complexity, both in terms of design and implementation. This increase is often significant enough to avoid using the KWT altogether and rather run the SPRT, on the grounds of its better nominal performance, or the FSST, on the grounds of it being simpler to implement and having an even smaller worst case sample size. In fact, the FSST can be considered ``maximally robust'', in the sense that its sample size is entirely independent not only of unknown parameters but of the sample distribution in general. 

In light of this discussion, the question addressed in this paper can be put as follows: Given a bound on the maximum expected sample size that is required to hold for all possible distributions, can one do better than a fixed-sample-size test? On the one hand, this question is interesting from a robustness perspective since it considers the expected sample size under a maximally large uncertainty set. On the other hand, it also touches upon the question of how fundamentally sequential tests and fixed-sample-size tests differ. More precisely, which additional constraints can be introduced in the design of a sequential tests before it ``collapses'' to a fixed-sample-size test?

\section{The Nonparametric Kiefer--Weiss Problem} 
\label{sec:nonpa_kiwei}                           

The hypothesis testing problem investigated in the paper is given by
\begin{equation}
  \min_{\pi \in \Pi} \; \sup_{\Pbb_0 \in \mathbb{M}_\mu} \; \gamma(\pi, \mathbb{P}_0) + \lambda_1 \alpha_1(\pi, \mathbb{P}_1) + \lambda_2 \alpha_2(\pi, \mathbb{P}_2).
  \label{eq:nonpa_kiwei}
\end{equation}
with $\lambda = (\lambda_1, \lambda_2) > 0$. As mentioned before, the problem in \eqref{eq:nonpa_kiwei} is a special case of the problem studied in \cite{Fauss2020_aos}, where the minimax formulation also includes the error probabilities, and the uncertainty sets are allowed to be arbitrary. Consequently, the results presented here heavily rely on the results in \cite{Fauss2020_aos}, and the reader will be referred to this reference for some of the proofs.

Before entering a more detailed discussion of how the solution of \eqref{eq:nonpa_kiwei} can be obtained from the results in \cite{Fauss2020_aos}, a necessary optimality condition is stated in the next Theorem. This condition is independent of the results in \cite{Fauss2020_aos} and provides some first insights into properties of the optimal test and the least favorable distributions.

\begin{theorem}[Necessary Optimality Conditions]
  In order for a pair $(\pi^*, \Qbb_0)$ to be optimal in the sense of \eqref{eq:nonpa_kiwei}, it needs to hold that
  \begin{align}
    \Ebb_{\pi^*} \bigl[ \tau_{\pi^*}(\bm{x}) \bigr] &\leq c \quad \forall \bm{x} \in \Xcal^\Nbb, \label{eq:inequality_condition} \\
    \Ebb_{\pi^*}\bigl[ \tau_{\pi^*}(\bm{x}) \bigr] &= c \quad    \mathbb{Q}_0 \text{-a.e.} \label{eq:equality_condition},
  \end{align}
  for some $c < \min\{\lambda_1, \lambda_2\}$, where $\Ebb_{\pi}$ denotes the expected value taken w.r.t.\ the randomization of the policy $\pi$.
  \label{th:optimality_necessary}
\end{theorem}

Theorem~\ref{th:optimality_necessary} is proven in Appendix~\ref{apx:optimality_necessary}. The inequality in \eqref{eq:inequality_condition} simply states that the worst-case expected sample size of the optimal test is bounded, and the equality in \eqref{eq:equality_condition} states that the least favorable distribution is concentrated on a set of sequences that attain this bound. This in turn implies that the expected sample size under $\Qbb_0$ is independent of the realization of $\bm{X}$. This equalization property is characteristic of minimax procedures in general \cite{VerduPoor1984}. 

It is not hard to show that any feasible FSST with sample size $N < \min\{\lambda_1, \lambda_2\}$ satisfies the conditions in Theorem~\ref{th:optimality_necessary} with $c = N$. Moreover, every distribution is least favorable w.r.t.\ the FSST in the sense of Theorem~\ref{th:optimality_necessary}. However, it is important to note that \eqref{eq:equality_condition} does not imply that the sample size of the optimal test is constant. First, $\Qbb_0$ might not be supported on the entire sample space, so that there can be sequences that lead to sample sizes smaller than $c$. Second, the same realization $\bm{x}$ can lead to different sample sizes, depending on the outcomes of the randomization; recall the definition in \eqref{eq:stopping_time}. 

For the SPRT and the standard KWT the randomization of the stopping rule is usually of little interest since both can be implemented using deterministic stopping rules.\footnote{%
  More precisely, the randomization can be arbitrary unless there is a positive probability for the test statistic to hit the thresholds exactly. However, even in this case, all randomization rules lead to the same weighted sum cost, meaning that the overall cost cannot be reduced, but only ``shifted'' between the expected sample size and the error probabilities; compare the discussion in \cite{Novikov2009, Fauss2015, Fauss2020_aos}.
}
For the nonparametric KWT (NP-KWT), however, it will become clear in the next section that randomized stopping rules are not critical. 

In order to obtain sufficient optimality conditions, the problem in \eqref{eq:nonpa_kiwei} is embedded into the more general problem formulation in \cite{Fauss2020_aos} by setting the number of hypotheses to two and choosing the uncertainty sets as
\begin{align}
  \Pcal_0 &= \Mcal_\mu, & \Pcal_1 &= \{P_1\}, & \Pcal_2 &= \{P_2\}.
  \label{eq:uncertainty_sets}
\end{align}
It then follows from Theorems~7 and 8 in \cite{Fauss2020_aos} that the optimal test in the sense of \eqref{eq:nonpa_kiwei} is characterized by the cost function $\rho_\lambda \colon \Rbb_+^3 \to \Rbb_+$, which is the unique solution of the equation
\begin{align}
  \rho_\lambda(\bm{z}) = \min \bigl\{\, & g_\lambda(\bm{z}) \,,\, z_0 + d_\lambda(\bm{z}) \,\bigr\},
  \label{eq:rho_implicit}
\end{align}
where $g_\lambda \colon \Rbb_+^3 \to \Rbb_+$ is given by
\begin{equation}
  g_\lambda(\bm{z}) \coloneqq \min\{\, \lambda_1 z_1 \,,\, \lambda_2 z_2 \,\},
\end{equation}
and $d_\lambda \colon \Rbb_+^3 \to \Rbb_+$ depends on $\rho_\lambda$ via
\begin{equation}
  d_\lambda(\bm{z}) = \sup_{P_0 \in \Mcal_\mu} \, \int_{\Xcal} \rho_\lambda \bigl( \bm{z} \bm{p}(x) \bigr) \, \mu(\mathrm{d} x).
  \label{eq:def_d}
\end{equation}
Here $\bm{p}$ is shorthand for the vector $(p_0, p_1, p_2)$. Also, note that $g_\lambda$ is independent of $z_0$. Depending on the context, both notations $g_\lambda(\bm{z})$ and $g_\lambda(z_1, z_2)$ are used in what follows. 

Alternatively, the function $\rho_\lambda$ can be defined as the limit
\begin{equation}
  \rho_\lambda = \lim_{n \to \infty} \; \rho_\lambda^{(n)},
\end{equation}
where 
\begin{align}
  \rho_\lambda^{(n)}(\bm{z}) &= \min \bigl\{\, g_\lambda(\bm{z}) \,,\, z_0 + d_\lambda^{(n)}(\bm{z}) \,\bigr\},
  \label{eq:rho_recursive} \\[1ex]
  d_\lambda^{(n)}(\bm{z}) &= \sup_{P_0 \in \Mcal_\mu} \, \int_{\Xcal} \rho_\lambda^{(n-1)} \bigl( \bm{z} \bm{p}(x) \bigr) \, \mu(\mathrm{d} x),
  \label{eq:def_dn}
\end{align}
and $\rho_\lambda^{(0)} = g_\lambda$. Here, $\rho_\lambda^{(n)}$ is the cost of an optimal test with finite horizon $n$, meaning that the test uses at most $n$ samples.\footnote{%
  The horizon of sequential test can be but does not have to be identical to its maximum sample size; the latter is a property of the test itself, the former is a design parameter. For example, the designer might choose a horizon of \num{50} samples, when the optimal test is in fact truncated after \num{30} samples; compare the KWT in Figure~\ref{fig:thresholds}.
}  
The optimal infinite-horizon test is then obtained by letting $n$ go to infinity. See \cite{Novikov2009, Fauss2020_aos} for a proof and a more detailed discussion.

In what follows, it is assumed that the suprema on the right-hand sides of \eqref{eq:def_d} and \eqref{eq:def_dn} are attained. This assumption is guaranteed to be satisfied for discrete sample spaces \cite{Fauss2018_tsp}. In light of the properties of $\rho_\lambda$ shown later in this section, we conjecture that it is also satisfied for a larger class of sample spaces. However, a formal characterization is bound to be technical, compare \cite{Rockafellar1967, Rockafellar1968, Papageorgiou1997}, and is beyond the scope of this paper. Also, note that if the supremum is not attained, $\rho_\lambda$ is still well-defined, but the corresponding optimal test is not, since its stopping rule depends on the least favorable distribution.

The importance of $\rho_\lambda$ lies in the fact that it provides sufficient and, in a sense, constructive conditions for minimax optimality. Before going into details of the latter, some useful properties of $\rho_\lambda$ and $\rho_\lambda^{(n)}$ are given in the next theorem.

\begin{theorem}
  [Properties of $\rho_\lambda$ and $\rho_\lambda^{(n)}$]
  Let $\bigl( \rho_\lambda^{(n)} \bigr)_{n \in \Nbb_0}$ denote the sequence of functions defined in \eqref{eq:rho_recursive} with limit $\rho_\lambda$. For all $n \in \Nbb_0$ it holds that
  \begin{enumerate}
      \item $\rho_\lambda^{(n)}$ and $\rho_\lambda$ are nondecreasing and concave.
      \item the right derivatives of $\rho_\lambda^{(n)}$ and $\rho_\lambda$ w.r.t.\ $z_0$ take values in $\Nbb_0 \cup \{ \infty \}$, more precisely,
            \begin{align}
              \partial_{z_0}^+ \rho_\lambda^{(n)}(\bm{z}) &\in \{0, \ldots, n\}, \\
              \partial_{z_0}^+ \rho_\lambda(\bm{z}) &\in \Nbb_0 \cup \{ \infty \}
            \end{align}
            for all $\bm{z} \in \mathbb{R}_+^3$. Moreover, for all $(z_1, z_2) \in \mathbb{R}_+^2$ it holds that
            \begin{align}
              \rho_\lambda(0, z_1, z_2) &= 0
              \label{eq:zero_sample_cost_rho}, \\
              \partial_{z_0}^+ \rho_\lambda(0, z_1, z_2) &\in \{0, \infty\}
              \label{eq:zero_sample_cost_drho}.
            \end{align}
  \end{enumerate}
  \label{th:rho_properties}
\end{theorem}

Theorem~\ref{th:rho_properties} is proven in Appendix~\ref{app:proof_th_2}. In words, it states that for any given $(z_1, z_2) \in \mathbb{R}_+^2$ the function $\rho_\lambda^{(n)}(\bullet, z_1, z_2)$ is piecewise linear with all segments admitting an integer valued slope of at most $n$. As a consequence, the partial derivative of $\rho_\lambda^{(n)}$ w.r.t.\ $z_0$ is a nonincreasing step function. More precisely, the partial differential at any point is either unique and integer valued, or it spans the ``size of the step'' at points where the right partial derivative is discontinuous. The limit $\rho_\lambda(\bullet, z_1, z_2)$ is piecewise linear on a possibly infinite partition of $\mathbb{R}_+$ with all segments admitting an integer valued but unbounded slope. 

A special case arises at $z_0 = 0$, where $\rho_\lambda$ is zero and its partial derivative w.r.t.\ $z_k$ can only take two values, either zero or infinity. Intuitively speaking, this effect occurs since $z_0$ in \eqref{eq:rho_implicit} and \eqref{eq:rho_recursive} corresponds to the cost for taking another sample; compare \cite{Novikov2009, Novikov2009_multiple_hypotheses, Fauss2015}. For $z_0 = 0$, this cost becomes zero so that the test can take additional samples ``for free'', thus reducing the error probabilities to zero. This, in turn, leads to an unbounded expected sample size ($\partial_{z_0}^+ \rho_\lambda = \infty$), unless the test has already stopped ($\partial_{z_0}^+ \rho_\lambda = 0$).

From an operational perspective, the vector $\bm{z}$ is the sufficient statistic of the optimal test. More precisely, it is shown in the next theorem that an optimal policy $\pi^*$ can be defined as a function of $\bm{z}^n$, $n \geq 0$, where
\begin{align}
  z_0^n = z_0^n(\bm{x}) &= z_0^{n-1} \frac{\mathrm{d} \Qbb_0}{\mathrm{d} \mu^\Nbb}(x_n | \bm{x}_{n-1}) \label{eq:z0} \\
  z_1^n = z_1^n(\bm{x}) &= z_1^{n-1} p_1(x_n) \label{eq:z1} \\
  z_2^n = z_2^n(\bm{x}) &= z_2^{n-1} p_2(x_n) \label{eq:z2}
\end{align}
and the initial value is denoted by $\bm{z}^0$. When running the test, $\bm{z}^0$ needs to be chosen as $\bm{z}^0 = (1, 1, 1)$ such that the elements of $\bm{z}^n$ correspond to the likelihood of $\bm{x}_n$ under $\Qbb_0$, $\Pbb_1$, and $\Pbb_2$, respectively. However, in order to characterize and analyze the test, it is useful to look at the expected sample size as a function of $\bm{z}^0$. To this end, let $\gamma_{\pi, \Pbb} \colon \Rbb_+^3 \to \Rbb_+$ be defined as
\begin{equation}
  \gamma_{\pi, \Pbb}(\bm{z}) = \Ebb_{\pi, \Pbb}\bigl[ \tau_\pi(\bm{X}) \mid \bm{z}^0 = \bm{z} \bigr],
  \label{eq:gamma_z}
\end{equation}
that is, $\gamma_{\pi, \Pbb}(\bm{z})$ denotes the expected sample size of a test with policy $\pi$ and initial test statistic $\bm{z}^0 = \bm{z}$ under distribution $\Pbb$. A sufficient condtion for $(\pi^*, \Qbb_0)$ to be minimax optimal can now be stated in terms of $\gamma_{\pi, \Pbb}$.

\begin{theorem}[Sufficient Optimality Conditions] 
  Let
  \begin{equation}
    \delta^*(\bm{z}) \begin{dcases}
                        = 1, & \lambda_1 z_1 > \lambda_2 z_2 \\
                        \in [0, 1], & \lambda_1 z_1 = \lambda_2 z_2 \\
                        = 0, & \lambda_1 z_1 < \lambda_2 z_2 \\
                      \end{dcases},
    \label{eq:delta_opt}
  \end{equation}
  \begin{equation}
    \psi^*(\bm{z}) \begin{dcases}
                        = 1, & g_\lambda(\bm{z}) < z_0 + d_\lambda(\bm{z}) \\
                        \in [0, 1], & g_\lambda(\bm{z}) = z_0 + d_\lambda(\bm{z}) \\
                        = 0, & g_\lambda(\bm{z}) > z_0 + d_\lambda(\bm{z}) \\
                      \end{dcases},
    \label{eq:psi_opt}
  \end{equation}
  and let $\bm{z}^n$ be as in \eqref{eq:z0}--\eqref{eq:z2} with $\bm{z}^0 = (1, 1, 1)$. In order for the pair $(\pi^*, \Qbb_0)$ to be optimal in the sense of \eqref{eq:nonpa_kiwei} it suffices that $\delta_n^*(\bm{x}) = \delta^*(\bm{z}^n)$, $\psi_n^*(\bm{x}) = \psi^*(\bm{z}^n)$, that
  \begin{equation}
  	\Qbb_0 = \prod_{n \in \Nbb_0} Q_{\bm{z}^{n}},
  	\label{eq:lfd_product}
  \end{equation}
  and that $Q_{\bm{z}}$ and $\psi^*(\bm{z})$ are such that
  \begin{align}
    c(\bm{z}) &\geq \gamma_{\pi^*, \Qbb_0} \bigl( z_0 q_{\bm{z}}(x), z_1 p_1(x), z_2 p_2(x) \bigr) \quad \forall x \in \Xcal 
    \label{eq:q_inequality} \\
    c(\bm{z}) &= \gamma_{\pi^*, \Qbb_0} \bigl( z_0 q_{\bm{z}}(x), z_1 p_1(x), z_2 p_2(x) \bigr) \quad Q_{\bm{z}}\text{-a.e.}  
    \label{eq:q_equality}
  \end{align}
  for some some function $c \colon \Rbb_+^3 \to \Nbb_0$.
  \label{th:optimality_sufficient}
\end{theorem}

Theorem~\ref{th:optimality_sufficient} is a specialized version of Theorem~8 in \cite{Fauss2020_aos} and is proven in Appendix~\ref{app:proof_th_3}. First, it states that, unsurprisingly, the optimal policy of the NP-KWT consists of a cost minimizing stopping rule and a likelihood-ratio-based decision rule. Second, it states that under the least favorable distribution, $\Qbb_0$, the process $\bm{X}$ is a Markov chain with sufficient statistic $\bm{z}^n$. The interesting part of Theorem~\ref{th:optimality_sufficient} is how the conditional distributions of this Markov chain are characterized, namely, in terms of equalization properties similar to those in Theorem~\ref{th:optimality_necessary}. However, while the conditions in \eqref{eq:inequality_condition} and \eqref{eq:equality_condition} merely require the overall sample size to be bounded and constant $\Qbb_0$-a.e., the conditions in \eqref{eq:q_inequality} and \eqref{eq:q_equality} take the dynamics of the sequential test into account: The optimal stopping rule and the least favorable distribution need to be such that in any given state $\bm{z}^n = \bm{z}$ the test can only transition to states with \emph{identical expected remaining sample size} $c(\bm{z})$.

For the finite-horizon case, the optimal stopping and decision rules become dependent on $n$ so that $\rho_\lambda$ and $d_\lambda$ need to be replaced by $\rho_\lambda^{(n)}$ and $d_\lambda^{(n)}$, respectively. Moreover, a corresponding sequence of functions $\gamma_{\pi^*, \Qbb_0}^{(n)}$ needs to be defined in order to state the finite-horizon versions of \eqref{eq:q_inequality} and \eqref{eq:q_equality}. Other than that, the proof carries over unchanged. The details provide little additional insight and are omitted for brevity. 

At first glance, it might seem that the equalization property in Theorem~\ref{th:optimality_sufficient} enforces a fixed-sample-size test: If the expected sample size in the initial state $\bm{z}^0$ is $N$, then the remaining expected sample size in ($\Qbb_0$-almost) every state $\bm{z}^1$ needs to be $N-1$, and so on, leaving seemingly no room for randomness in sample size. However, as will become clear soon, the NP-KWT ``creates'' randomness by exploiting the option to randomize its stopping rule. In other words, while the sample size of the SPRT or the KWT is random because different sequences correspond to different stopping times, the sample size of the NP-KWT is random because the stopping time of individual sequences is random. This difference will become more tangible in Section~\ref{sec:example}, where it is illustrated with a numerical example. Before that, some useful general properties of the NP-KWT are discussed in the next section.

\section{Properties of the Nonparametric Kiefer--Weiss Test} 
\label{sec:nonpa_kiwei_properties}                           

The aim of this section is to provide a higher-level perspective on the NP-KWT and, where possible, to given an intuitive understanding of its properties. Note that, throughout this section, when speaking of an NP-KWT, what is really meant is an NP-KWT \emph{with policies according to Theorem~\ref{th:optimality_sufficient}}.

First, the optimal stopping rule characterized in Theorem~\ref{th:optimality_sufficient} can be stated more explicitly in terms of a threshold test, which is in closer analogy to regular sequential tests.

\begin{corollary}
  The optimal stopping rule $\psi^*$ in Theorem~\ref{th:optimality_sufficient} is of the form 
  \begin{equation}
    \psi^*(\bm{z}) \begin{dcases}
                    = 0, & z_0 < z_0^* \\
                    \in [0, 1], & z_0 = z_0^* \\
                    = 1, & z_0 > z_0^*
                   \end{dcases}
    \label{eq:psi_opt_threshold}
  \end{equation}
  where $z_0^* = z_0^*(z_1, z_2)$ is unique and depends on $(z_1, z_2)$. 
  \label{cl:threshold_test}
\end{corollary}

Corollary~\ref{cl:threshold_test} is proven in Appendix~\ref{app:proof_cl_1}. Rewriting the stopping rule in this form allows for some insight into the underlying policy. On the one hand, the NP-KWT stops if the likelihood of the observed sequence under the least favorable distribution becomes sufficiently large. This corresponds to cases in which the true distribution is likely not to be any of the nominal distributions, hence there is no use in continuing the test. On the other hand, the test stops if the threshold becomes sufficiently small. The threshold value is in turn determined by the likelihood of the observed sequence under $\Pbb_1$ and $\Pbb_2$. That is, the threshold for stopping is increased or decreased, depending on the evidence for the hypotheses. Of course, when running the test, both effects are coupled and jointly influence the stopping time. However, it is useful to keep the distinction in mind between stopping because there is sufficient evidence for a hypothesis and stopping because continuing is likely to be futile. 

An interesting property of the least favorable distribution is stated in the next corollary.

\begin{corollary}
  Let $\Qcal_{\bm{z}}$ denote the set of distributions that are least favorable in the sense of Theorem~\ref{th:optimality_sufficient} for a given value of $\bm{z}$. For all $\bm{z} \in \Rbb_+^3$, one of the following two statements holds true:
  \begin{enumerate}
    \item $\Qcal_{\bm{z}} = \Mcal_\mu$.
  	\item $\supp(Q) \supset \supp(P_1) \cap \supp(P_2)$ for all $Q \in \Qcal_{\bm{z}}$.
  \end{enumerate}
  Moreover, the inclusion in the second statement is strict only for the last sample. That is, $X_n \sim Q$ and $\supp(Q) \varsupsetneq \supp(P_1) \cap \supp(P_2)$ implies $\tau_{\pi}(\bm{x}) = n$ for all $x_n \in \Xcal$.
  \label{cl:lfd_support}
\end{corollary}

Corollary~\ref{cl:lfd_support} is proven in Appendix~\ref{app:proof_cl_2}. It states that every conditionally least favorable distribution is supported on the intersection of the supports of $P_1$ and $P_2$, with the exception of trivial cases in which every distribution is conditionally least favorable. In other words, events that occur with probability zero under either $P_1$ or $P_2$ can only occur under $Q$ if the underlying test has already stopped or is guaranteed to stop after the next sample. The latter is a corner case that will be discussed in more detail later in this section. 

For the common case where $\supp(P_1) = \supp(P_2) = \Xcal$, Corollary~\ref{cl:lfd_support} implies that $\Qbb_0$ can always be chosen such that it is supported on the entire sample space $\Xcal^\Nbb$. This property might seem counterintuitive since one would expect the least favorable distribution to exclude sequences that are highly indicative of either hypothesis. However, it is clear from Corollary~\ref{cl:threshold_test} that the NP-KWT also stops if a sequence has a high likelihood under the least favorable distribution. Hence, loosely speaking, the more concentrated $\Qbb_0$ is, the easier it can be recognized.

Picking up the coin flipping example from the Prologue, the intuition that $\Qbb_0$ should be concentrated on a subset of nonindicative sequences corresponds to Sabri's first attempt at prolonging the test by following a deterministic pattern. In this case, $\Qbb_0$ is maximally concentrated, which in turn makes it trivial for Beenish to identify it. Hence, even from a purely intuitive perspective, $\Qbb_0$ needs to be supported on a set of sequences that is ``large enough'' to mimic the nominal distributions. Corollary~\ref{cl:lfd_support} makes this notion precise: the support is ``large enough'' if it contains \emph{all events that occur with positive probability under both hypotheses}. Of course, this does not imply that all sequences are equally likely. As will be shown in the next section, events that are highly indicative of $\Hcal_1$ or $\Hcal_2$ are indeed less likely to occur under $\Qbb_0$ than under the respective hypothesis.

Finally, note that Corollary~\ref{cl:lfd_support} does not hold for tests with finite horizon, that is, for tests whose behavior is governed by the sequence $\bigl( \rho_\lambda^{(n)} \bigr)_{n \geq 0}$ instead of the limit $\rho_\lambda$. The reason for this is that the partial derivative of $\rho_\lambda^{(n)}$ w.r.t.\ $z_0$, is bounded away from infinity. As a consequence, the expected remaining sample size is bounded, which makes it feasible for the optimal test to enter the respective states. However, since $\rho_\lambda^{(n)}$ converges to $\rho_\lambda$, the probability of an optimal test entering a state in which Corollary~\ref{cl:lfd_support} does not hold goes to zero for $n \to \infty$. 

The fact that the properties of the least favorable distributions depend on the horizon of the test points to another question of interest, namely, the question of whether or not the NP-KWT with infinite horizon is truncated. To clarify, a sequential test with policy $\pi$ is referred to as being truncated if there exists some $N \geq 0$ such that
\begin{equation}
  \esssup \, \tau_\pi(\bm{x}) \leq N \quad \text{for all} \ \bm{x} \in \Xcal^\Nbb,
  \label{eq:def_truncated}
\end{equation}
where the essential supremum is taken with respect to the distribution induced by the randomization of the stopping rule. In words, a test is truncated in the above sense if the probability of its sample size exceeding $N$ is zero for all possible sequences of observations.

While the SPRT is nontruncated, which causes its expected sample size to be sensitive to deviations from the nominal distributions, both the KWT and, trivially so, the FSST are truncated in the above sense. This suggests that the NP-KWT, being positioned between the KWT and the FSST in terms of its robustness properties, might also be truncated. Interestingly, it turns out that this is not necessarily the case. 

\begin{corollary}
  The NP-KWT in the sense of Theorem~\ref{th:optimality_sufficient} is in general neither truncated nor nontruncated. A sufficient condition for it to be truncated is that
  \begin{itemize}
    \item $\supp(P_1) \cap \supp(P_2)$ is a singleton.
  \end{itemize}
  A sufficient condition for it to be nontruncated is that
  \begin{itemize}
    \item $\supp(P_1) = \supp(P_2)$ and
    \item $\lambda_1 = \lambda_2 = \lambda$ large enough such that $\tau_{\pi^*}(\bm{x}) > 1$ for all $\bm{x} \in \Xcal^\Nbb$.
  \end{itemize}
  \label{cl:truncation}
\end{corollary}

Corollary~\ref{cl:truncation} follows from Theorem~\ref{th:rho_properties} and is proven in Appendix~\ref{app:proof_cl_3}. The respective conditions on the intersection of the supports can be seen as extreme cases. Qualitatively speaking, the NP-KWT is truncated if the support of the LFDs is too small, that is, if there are not enough degrees of freedom to choose the LFDs in a manner that prevents the test form stopping with certainty. Building on the arguments in Appendix~\ref{app:proof_cl_3}, it should be possible to make this notion more precise, but this exercise will not be attempted here.

The conditions on $\lambda$ in the second statement of the corollary can also be relaxed. In general, $\lambda$ needs to be such that it is possible for the test to reach a state in which the probability of the next sample overwriting the evidence collected so far is sufficiently small. For symmetric cost coefficients, $\lambda_1 = \lambda_2$, such a state can always be reached with the first sample, which leads to the conditions given above. For $\lambda_1 \neq \lambda_2$, the NP-KWT is also likely to be truncated, but it might require a larger minimum sample size, depending on symmetry properties of the likelihood ratio. In order to keep the corollary succinct, this case is not included. Again, a refinement should be possible based on the arguments outlined in Appendix~\ref{app:proof_cl_3}.

The fact that the NP-KWT can be nontruncated makes sense in light of the optimality conditions in Theorem~\ref{th:optimality_sufficient}. Namely, as discussed before, \eqref{eq:q_inequality} and \eqref{eq:q_equality} imply that at any instant $n$ the NP-KWT can only transition to states with identical expected remaining sample size, irrespective of the realization of the next sample $X_{n+1}$. This means that the decision to stop has to be made \emph{before} $X_{n+1}$ is observed since stopping for one value of $X_{n+1}$ implies stopping for ($\Qbb_0$-almost) \emph{every} value of $X_{n+1}$. This creates a situation in which the very objective of the test is to stop as early as possible, yet every stopping decision comes with the risk of not knowing the last sample. The NP-KWT solves this problem by \emph{probably} stopping after the next sample. In this way, it does not need to commit to stopping in cases where the next sample leads to a more ambiguous state, yet is bound to stop eventually after a sequence of indicative samples. 
As suggested by the title of the paper, the test does not cross a hard boundary, but softly fades out as the evidence increases.

\section{Example: Coin Flipping} 
\label{sec:example}              

In this section, the properties of the NP-KWT with policies according to Theorem~\ref{th:optimality_sufficient} are illustrated using the coin flipping example from the Prologue. More formally speaking, the aim is to design an NP-KWT for the success probability of a Bernoulli random variable. This example was chosen for several reasons. First, given that $\lambda$ is chosen appropriately, it satisfies the sufficient conditions in Corollary~\ref{cl:truncation} for the NP-KWT to be nontruncated, which is arguably the more interesting case. Second, it is a simple, tractable example. Although the least favorable distribution and the optimal policy in Theorem~\ref{th:optimality_sufficient} can in principle be calculated explicitly, doing so is non-trivial and requires an adequate way of representing $\rho_\lambda$ and $\rho_\lambda^{(n)}$. Third, the sequences of observations generated by flipping coins are easy to visualize using binary trees. This is important since there is no obvious graphical representation of the NP-KWT that is comparable to those in Figure~\ref{fig:thresholds} in terms of usefulness and clarity. Finally, the Bernoulli example can be considered to be the most basic scenario in the sense that every other test can be reduced to it by partitioning the sample space into two subsets.

Consider the hypotheses in \eqref{eq:hyp_illustration}, with corresponding PMFs $p_1$ and $p_2$ in \eqref{eq:nominal_densities}. In order to implement the NP-KWT, the function $\rho_\lambda$ in \eqref{eq:rho_implicit} needs to be known. However, calculating or approximating $\rho_\lambda$ is a non-trivial problem. In non-robust sequential detection, similar cost functions are guaranteed to be smooth, or can be represented as the minimum of two smooth functions, which means that they can be approximated using standard techniques, such as splines or basis functions \cite{Fauss2015, Fauss2018_tsp}. These techniques cannot be applied in a straightforward manner here. The piecewise linear nature of $\rho_\lambda$, including the exact locations of the changes in slope, is crucial and needs to be preserved in order for the corresponding test to be meaningful.

The numerical results presented here are based on a finite-horizon approximation of $\rho_\lambda$. That is, for some $N > 0$, the sequence $\bigl( \rho_\lambda^{(n)} \bigr)_{0 \leq n \leq N}$ is calculated exactly and recursively using \eqref{eq:rho_recursive}. This procedure is somewhat unsatisfactory, given that most of the properties stated in the previous section only hold in the limit. However, the tests obtained in this manner are strictly optimal for the given horizon and are guaranteed to converge to the infinite-horizon case.

The coin flipping example also allows for some additional simplifications in the numerical design. In particular, both likelihoods $z_1^n$ and $z_2^n$ can be reparametrized in terms of the number of successes or failures observed so far. In this way, $\rho_\lambda^{(n)}(z_0, z_1, z_2)$, with domain $\Rbb_+^3$, can be represented by a family of functions $\rho_\lambda^{(n,m)}(z_0)$ with domain $\Rbb_+$, where $m$ denotes the number of successes among the first $n$ observations. For more implementation details, the interested reader is referred to the repository \cite{github}, where Python code for the design of finite-horizon NP-KWTs for discrete distributions can be found.

\begin{figure}
  \centering
  \includegraphics{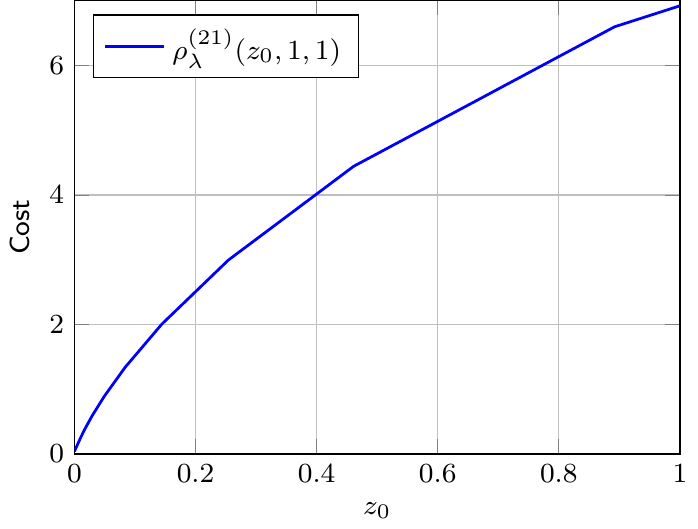} \hfill
  \includegraphics{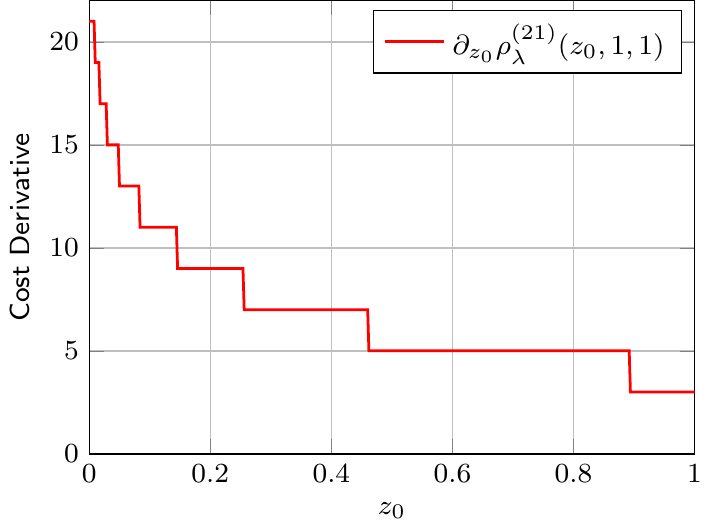}\\
  \caption{Example for an optimal cost function $\rho_\lambda^{(n)}$ and its partial derivative $\partial_{z_0} \rho_\lambda^{(n)}$ for the hypotheses in \eqref{eq:hyp_illustration}. Here, $n = 21$ and $\rho_\lambda^{(n)}$ has been calculated recursively according to \eqref{eq:rho_recursive}.}
  \label{fig:rho}
\end{figure}

For $n = 21$ and $\lambda_1 = \lambda_2 = 20$, an example of the function $\rho_\lambda^{(n)}$ and its partial derivative w.r.t.\ $z_0$ is shown in Figure~\ref{fig:rho}. Note that since $p_1, p_2 < 1$ in this example, $\rho_\lambda^{(n)}$ only needs to be evaluated on the interval $[0, 1]$. The piecewise linear form of $\rho_\lambda^{(n)}$ can be gauged from the left plot, but is more obvious from the partial derivative shown in the right plot. In accordance with Theorem~\ref{th:rho_properties}, the latter is non-increasing, only takes on integer values, and is bounded by $n = 21$, which is attained at $z_0 = 0$. The fact that it decreases with step-sizes of two is a side-effect of the underlying Bernoulli model, which leads to optimal tests admitting an odd expected sample size.

\begin{sidewaysfigure*}
  \centering
  \includegraphics[width=\textwidth]{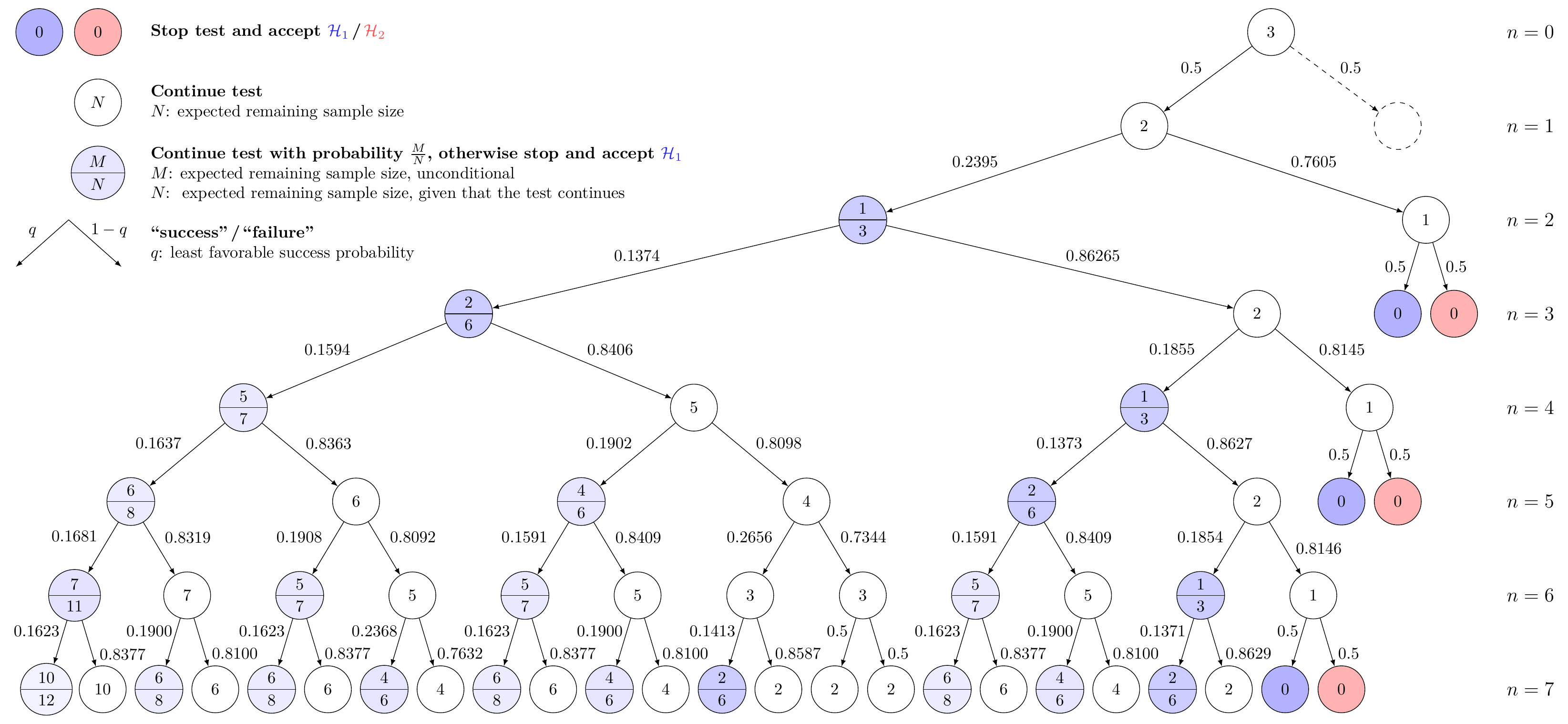}
  \caption{Optimal testing policy and least favorable distribution for the first seven samples of a nonparametric Kiefer--Weiss test for the hypotheses in \eqref{eq:hyp_illustration} with horizon $N = 21$ and maximum expected sample size of \num{3}. The test starts in the root node at the top and transitions to the next state based on whether a ``success'' (left arrow) or a ``failure'' (right arrow) is observed. The label(s) in each node contains information about the expected remaining sample size, given that the respective node has been reached. Three different cases need to be distinguished: White nodes correspond to states in which the test is guaranteed to continue. Here the label indicates the expected remaining sample size. Colored nodes with label ``$0$'' correspond to states in which the test is guaranteed to stop, blue indicating a decision for $\Hcal_1$, red for $\Hcal_2$. Shaded nodes with two labels correspond to states in which the test stops with a certain probability, which is reflected in the intensity of the shading. The top label corresponds to the expected remaining sample size when entering the state, that is, before the randomized stopping rule is evaluated. The lower label corresponds to the expected remaining sample size given that the test continues. Hence, the probability of continuing is given the ratio of both values, that is, the labels interpreted as a fraction. The numbers next to the arrows denote the probabilities of the respective observation under the least favorable distribution $\Qbb_0$. To give an example, if the first two observations are successes, the test enters the left node on level $n=2$. Here, the test continues with probability $\tfrac{1}{3}$ or stops with a decision for $\Hcal_1$ with probability $\tfrac{2}{3}$. Given that it continues, its expected remaining sample size is \num{3}, and the least favorable success probability of the next sample is $\approx 0.1374$. Note that since the decision tree is symmetric, only its left half is depicted, that is, only sequences starting with ``success''.}
  \label{fig:decision_tree}
\end{sidewaysfigure*}

The NP-KWT corresponding to the cost function in Figure~\ref{fig:rho} is illustrated by the binary tree in Figure~\ref{fig:decision_tree} on the next page; see the legend and the caption text for a detailed explanation of how to read this figure. Note that the expected sample size of the depicted test is $\Ebb_{\pi^*,\Qbb_0}\bigl[\, \tau_{\pi^*}(\bm{X}) \,\bigr] = 3$. This small expected sample size is intentional since it facilitates the graphical representation of the test.

Figure~\ref{fig:decision_tree} illustrates the most prominent properties of the NP-KWT. The six states in which the test is guaranteed to stop correspond to cases in which it is likely that the samples were drawn from $\Qbb_0$. Moreover, note that in the second-to-last states in theses cases the same number of successes and failures has been observed so that there is zero evidence for either hypothesis. In light of the discussion in the previous section and in Appendix~\ref{app:proof_cl_3}, this is an example of a scenario where a stopping decision can be made safely in the second-to-last state, since the evidence can only increase with the next sample.

For the majority of states, however, the ``fading out'' effect caused by the need to randomize the stopping rule in cases where the next sample can either increase of decrease the evidence can be seen. For example, after having observed two successes in a row, it seems reasonable to stop the test if another success is observed. However, this implies also stopping the test if a failure is observed, which leads to a much more ambiguous state. As discussed in the previous section, the NP-KWT works around this dilemma by stopping with probability $\tfrac{2}{3}$ if another success is observed and continuing with certainty if a failure is observed. The same rationale applies to all shaded nodes in the tree.  

By inspection, the expected sample size of the NP-KWT is independent of the true distribution, which is the minimax property required by both Theorem~\ref{th:optimality_necessary} and Theorem~\ref{th:optimality_sufficient}. For every node in Figure~\ref{fig:decision_tree}, its two child nodes have identical expected remaining sample sizes. However, nodes with random stopping rules can in fact increase this number. Their expected remaining sample size, conditioned on a decision to continue the test, can grow far beyond the overall expected sample size. In this example, it becomes largest (\num{12}) after having observed seven successes in a row. Again, the effect that highly indicative states correspond to large expected remaining sample sizes seems counterintuitive at first glance, but makes sense in light of the previous discussion.

Another noteworthy aspect of the NP-KWT is that, in contrast to the SPRT or the KWT, it admits a kind of memory in the form of $z_0^n$. For example, consider the sequences ``success-failure-success'' and ``success-success-failure''. Both sequences correspond to the same likelihood of $\Hcal_1$ being true, yet the former results in a definite decision for $\Hcal_1$, while the latter results in a test that will continue for at least two more samples. Also, note that this difference is ``the wrong way around'', in the sense that the sequence with stronger initial evidence, ``success-success'', leads to a larger sample size. Such a scenario cannot occur for the SPRT and the KWT. This difference is due to the fact that the SPRT and the KWT only stop if there is \emph{sufficient evidence} for either hypothesis, while the NP-KWT also stops if there is \emph{insufficient evidence} for either hypothesis.

Finally, Figure~\ref{fig:decision_tree} also shows the least favorable transition probabilities for each node. Interestingly, the least favorable success probabilities are not extreme values close to zero or one, but are rather similar to those under the hypotheses, with the lowest success probability being $\approx$\SI{13.7}{\percent} and the highest being $\approx$\SI{81.5}{\percent}. Again, the underlying rationale is that the least favorable distribution needs to be ``disguised'' as a nominal distribution in order to prevent the test from stopping early. Note that since the expected remaining sample size remains constant in all states, the least favorable distribution does not stand out as being particularly harmful, but is just as good or bad as any other distribution in terms of the expected sample size. Nevertheless, $\Qbb_0$ needs to be known since the stopping rule is a function of $z_0^n$. 

An open question is in how far the policy illustrated in Figure~\ref{fig:decision_tree} allows for insights into the infinite-horizon case. Although there seem to be patterns emerging in the sub-trees, there still are subtle differences between the branches. In general, it seems difficult to identify a rule that could be used to extrapolate the optimal policy to larger sample sizes. 

\begin{figure}
  \centering
  \includegraphics{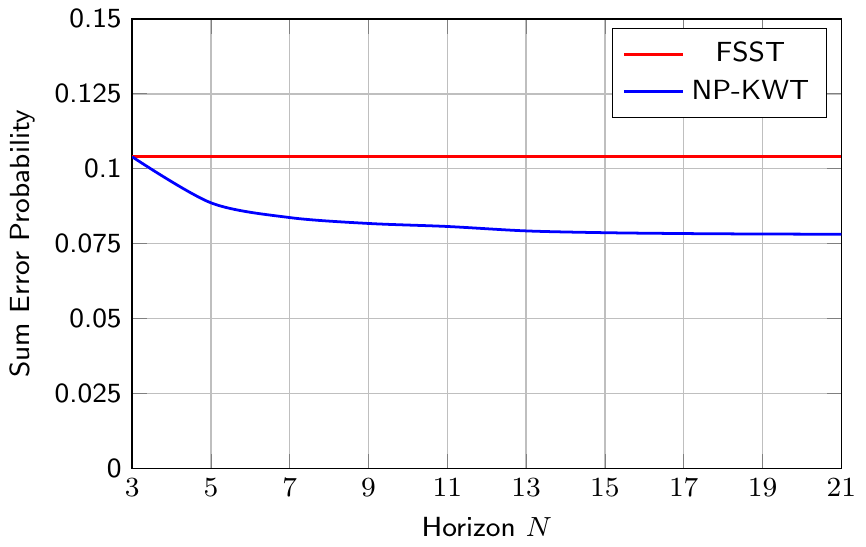}
  \caption{Comparison of the sum error probabilities of a fixed-sample-size test (FSST) and nonparametric Kiefer--Weiss test (NP-KWT) with horizon $N$. Both tests were designed to admit an expected sample size of three under all distributions $\Pbb \in \Mbb_\mu$.}
  \label{fig:NPKWT_vs_FSST}
\end{figure}

In order to put the NP-KWT into perspective, it can be compared to an FSST whose sample size is chosen to be identical to the expected sample size of the NP-KWT. In this case, both tests admit the same expected sample size under all distributions, so that the comparison is fair in this regard. In this example, the sum error probability, $\alpha_1(\pi, \mathbb{P}_1) + \alpha_2(\pi, \mathbb{P}_2)$, is used as a performance metric for the comparison. In Figure~\ref{fig:NPKWT_vs_FSST}, the sum error probability of both the FSST and the finite-horizon NP-KWT is plotted as a function of the horizon $N$. Note that all tests where designed such that they admit an expected sample size of three. It can be seen that with increasing horizons the NP-KWT outperforms the FSST by a small, but notable margin. The test illustrated in Figure~\ref{fig:decision_tree}, with $N = 21$, reduces the sum error probability from \SI{10.4}{\percent} for the FSST to $\approx$\SI{7.82}{\percent}, which is a reduction by almost \SI{25}{\percent}. This number will be slightly higher for the infinite-horizon NP-KWT, however, as can be seen in Figure~\ref{fig:NPKWT_vs_FSST}, the sum error probability appears to be close to convergence for $N = 21$ already. Finally, note that the probabilities shown in Figure~\ref{fig:NPKWT_vs_FSST} were calculated, not simulated, meaning that they are precise up to rounding errors.

\section{Conclusion and Outlook} 
\label{sec:conclusion}           

In this paper, a nonparametric variety of the Kiefer--Weiss problem has been formulated and two optimality conditions have been given, one necessary and one sufficient. Based on the latter, properties of a test that solves the nonparametric Kiefer--Weiss problem and the corresponding least favorable distributions have been derived. In order to illustrate the optimal testing policy, a nonparametric Kiefer--Weiss test for the success probability of a Bernoulli random variable has been implemented for a finite horizon, and has been shown to admit smaller error probabilities than a fixed-sample-size test with identical expected sample size. Moreover, the example has revealed that the nonparametric Kiefer--Weiss test heavily relies on randomized stopping rules, with deterministic stopping decisions being limited to sequences that are highly unlikely to occur under the hypothesized distributions.

Although the results in this paper allow for some insights into the working principles of the NP-KWT, they still leave a number of important questions unanswered.

\begin{itemize}
  \item In how far are the sufficient optimality conditions in Theorem~\ref{th:optimality_sufficient} also necessary? For example, the policies characterized by the conditions in Theorem~\ref{th:optimality_sufficient} are time-homogeneous, meaning that the decision and stopping rules are allowed to depend on the test statistic, $\bm{z}^n$, but not the sample number, $n$. Do optimal time-heterogeneous policies exist? This leads to the closely related questoin of whether the NP-KWT and the corresponding LFDs are unique. The latter can be shown to be the case for the Bernoulli example, but, given that $\rho_\lambda$ is not strictly convex, it is not clear whether this property holds in general.
  
  \item How can an infinite-horizon NP-KWT be designed and implemented? Clearly, a purely numerical design based on $\rho_\lambda$ is not going to be feasible. Hence, an alternative approach to the test design is likely to be necessary that does not require knowledge of $\rho_\lambda$. Such approaches exist for the SPRT \cite{Wald1947} and the KWT \cite{Robbins1972, Tartakovsky2014}, however, for the NP-KWT the design problem is significantly harder since it cannot be reduced to determining upper and lower thresholds for a scalar test statistic. In general, it would be useful to be able to characterize states in which the NP-KWT continues/stops/might stop without having to go through the entire design process.
\end{itemize}

Finally, we would like to conclude the paper with a conjecture. Clearly, the basis for the performance improvement of the NP-KWT over the FSST in the coin flipping example is the additional degree of freedom introduced by only requiring it to admit a constant \emph{expected} sample size, which makes it possible to occasionally observe long sequences. The increase in confidence gained from these sequences outweighs the slight decease in confidence incurred by sometimes having to stop the test short. This observation motivates the following conjecture.

\begin{conjecture}
  A fixed-sample-size test of length $N = c$, with $c$ as in Theorem~\ref{th:optimality_necessary},  solves \eqref{eq:nonpa_kiwei} if and only if
  \begin{itemize}
    \item $\supp(P_1) = \supp(P_2)$ and
    \item the stopping rule is constrained to be deterministic.
  \end{itemize}
\end{conjecture}

The conjecture is based on the following line of arguments: For deterministic testing policies, the necessary optimality conditions in Theorem~\ref{th:optimality_necessary} become
\begin{align}
  \tau_{\pi^*}(\bm{x}) &\leq c \quad \forall \bm{x} \in \Xcal^\Nbb, \\
  \tau_{\pi^*}(\bm{x}) &= c \quad    \mathbb{Q}_0 \text{-a.e.}
\end{align}
That is, $\Qbb_0$-almost all sequences need to admit the same \emph{deterministic} sample size. According to Corollary~\ref{cl:lfd_support}, $\Qbb_0$ is supported on $\supp(P_1) = \supp(P_2)$. Hence, the optimal test uses exactly $c$ samples with probability one under $\Qbb_0, \Pbb_1$ and $\Pbb_2$--- the fixed-sample-size test clearly admits this property.

The crucial step when following this sketch is to show that the support of the LFDs is independent of the randomization of the stopping rule. This should be the case since the LFDs are defined as maximizers of the function $\rho_\lambda$, which itself is independent of the randomization of the testing policy. However, a formal proof is beyond the scope of this paper.

\begin{appendix}

\section{Proof of Theorem~\ref{th:optimality_necessary}} 
\label{apx:optimality_necessary}                         

Both statements in Theorem~\ref{th:optimality_necessary} can be proven by contradiction. Assume that $(\pi^*, \Qbb_0)$ is minimax optimal. The objective function in \eqref{eq:nonpa_kiwei} is upper bounded by $\min\{\lambda_1, \lambda_2\}$, and the bound is attained by the sub-optimal policy
\begin{align}
  \psi_0 &= 1, & \delta_0 &= \begin{cases}
                              0, & \lambda_1 \leq \lambda_2 \\
                              1, & \lambda_1 > \lambda_2
                             \end{cases}.
\end{align}
Hence, there exists some $c < \min\{\lambda_1, \lambda_2\}$ such that
\begin{equation}
  \Ebb_{\pi^*, \Qbb_0}[\, \tau_{\pi^*}(\bm{X}) \,] \leq c.
  \label{eq:lfd_bounded}
\end{equation}
and
\begin{equation}
  \Ebb_{\pi^*}[\, \tau_{\pi^*}(\bm{x}') \,] = c.
  \label{eq:c_attained}
\end{equation}
for at least one sequence $x' \in \Xcal^\Nbb$. By definition of $\Qbb_0$, it holds that
\begin{equation}
  \Ebb_{\pi^*, \Qbb_0}[\, \tau_{\pi^*}(\bm{X}) \,] \geq \Ebb_{\pi^*, \Pbb}[\, \tau_{\pi^*}(\bm{X}) \,]
  \label{eq:lfd_condition}
\end{equation}
for all $\Pbb \in \Mbb^\Nbb$. Now, assume that a sequence $\tilde{\bm{x}}$ exists such that
\begin{equation}
  \Ebb_{\pi^*}[\, \tau_{\pi^*}(\tilde{\bm{x}}) \,] > c
\end{equation}
and let $\bm{1}_x$ denote the distribution with a single point mass of $1$ at $\bm{x}$. According to \eqref{eq:lfd_condition}, it holds that
\begin{align}
  \Ebb_{\pi^*, \Qbb_0}[\, \tau_{\pi^*}(\bm{X}) \,] &\geq \Ebb_{\pi^*, \bm{1}_{\tilde{x}}}[\, \tau_{\pi^*}(\bm{X}) \,] \\
  &= \Ebb_{\pi^*}[\, \tau_{\pi^*}(\tilde{\bm{x}}) \,]  \\
  &> c, \label{eq:lfd_contradiction}
\end{align}
which contradicts \eqref{eq:lfd_bounded} and hence $\pi^*$ being minimax optimal. This proves the fist statement.

A similar argument can be used to show the second statement. Let $\bm{x}'$ be as in \eqref{eq:c_attained}. It then holds that
\begin{align}
  \Ebb_{\pi^*, \Qbb_0}[\, \tau_{\pi^*}(\bm{X}) \,] &\geq \Ebb_{\pi^*, \bm{1}_{x'}}[\, \tau_{\pi^*}(\bm{X}) \,] \\
  &= \Ebb_{\pi^*}[\, \tau_{\pi^*}(\bm{x}') \,]  \\
  &= c.
\end{align}
In combination with \eqref{eq:lfd_bounded}, this implies the second statement in Theorem~\ref{th:optimality_necessary} and completes the proof.

\section{Proof of Theorem~\ref{th:rho_properties}} 
\label{app:proof_th_2}                             

The first statement in Theorem~\ref{th:rho_properties} is shown in \cite[Lemma 1]{Fauss2020_aos}. The second statement can be shown via induction. Assume that for some $n \geq 0$ it holds that $\rho_\lambda^{(n)}$ satisfies the properties stated in Theorem~\ref{th:rho_properties}. From \eqref{eq:rho_recursive} it follows that the partial differential of $\rho_\lambda^{(n+1)}$ w.r.t.\ $z_0$ evaluated at $\bm{z}$ is given by 
\begin{equation}
  \partial_{z_0}^+ \rho_\lambda^{(n+1)}(\bm{z}) = \begin{dcases}
                                                  0, & z_0 \geq z_0^* \\
                                                  1 + \partial_{z_0}^+ d_\lambda^{(n)}(\bm{z}), & z_0 < z_0^*  
                                                \end{dcases}
  \label{eq:derivative_rho_n}
\end{equation}
where $z_0^* = 0$ if $g_\lambda(\bm{z}) < z_0 + d_\lambda^{(n)}(\bm{z})$ for all $z_0 \geq 0$, and $z_0^*$ is the unique solution of
\begin{equation}
  g_\lambda(\bm{z}) = z_0 + d_\lambda^{(n)}(\bm{z})
  \label{eq:def_z0_star}
\end{equation}
otherwise. Uniqueness of $z_0^*$ follows from the fact that $g_\lambda$ is independent of $z_0$ and that the right-hand side of \eqref{eq:def_z0_star} is continuous, unbounded, and increasing in $z_0$. 

Recall that $d_\lambda^{(n)}$ is defined as
\begin{align}
  d_\lambda^{(n)}(\bm{z}) &= \sup_{P_0 \in \Mcal_\mu} \; \int_\Xcal \rho_\lambda^{(n)}(\bm{z} \bm{p}(x)) \, \mu(\mathrm{d}x) \\
  &= \int_\Xcal \rho_\lambda^{(n)}(\bm{z} \bm{q}(x)) \, \mu(\mathrm{d}x)
  \label{eq:d_recall}
\end{align}
where $\bm{q} = (q_0, p_1, p_2)$ is introduced for the sake of a more compact notation. A necessary and sufficient condition for a distribution $Q_0$ to solve \eqref{eq:d_recall}, is that there exists a constant $c$ and a partial derivative $r \in \partial_{z_0} \rho_\lambda^{(n)}$ such that 
\begin{align}
    r(\bm{z} \bm{q}(x)) &\leq c \quad \forall x \in \Xcal
    \label{eq:opt_condition_inequality} \\
    r(\bm{z} \bm{q}(x)) &= c \quad Q \text{-a.e.}
    \label{eq:opt_condition_equality}
\end{align}
These optimality conditions can be shown in analogy to the proof of Theorem~\ref{th:optimality_necessary}; also compare Theorem~1 in \cite{Fauss2018_tsp}. Let the set of functions $r$ that satisfy \eqref{eq:opt_condition_inequality} and \eqref{eq:opt_condition_equality} for a given $c$ be denoted by $\mathcal{R}_c(\bm{z})$, that is,
\begin{equation}
  \mathcal{R}_c^{(n)}(\bm{z}) \coloneqq \left\{ r \in \partial_{z_0} \rho_\lambda^{(n)} : \exists Q \in \Mcal_\mu : \eqref{eq:opt_condition_inequality} \ \text{and} \ \eqref{eq:opt_condition_equality} \ \text{hold} \right\}.
  \label{eq:R_set}
\end{equation}
Further, let $\mathcal{C}^{(n)}(\bm{z})$ denote the set of nonnegative scalars $c$ for which $\mathcal{R}_c^{(n)}(\bm{z})$ is non-empty:
\begin{equation}
  \mathcal{C}^{(n)}(\bm{z}) \coloneqq \left\{ c \in \Rbb_+ : \mathcal{R}_c^{(n)}(\bm{z}) \neq \emptyset \right\}.
\end{equation}
According to Lemma~3 in \cite{Fauss2020_aos}, differentiation and integration on the right-hand side of \eqref{eq:d_recall} can be interchanged. Hence, it holds that
\begin{align}
  \partial_{z_0}^+ d_\lambda^{(n)}(\bm{z}) &= \partial_{z_0}^+ \int_\Xcal \rho_\lambda^{(n)}(\bm{z} \bm{q}(x)) \, \mu(\mathrm{d}x) \\
  &= \min \; \left\{ \partial_{z_0} \int_\Xcal \rho_\lambda^{(n)}(\bm{z} \bm{q}(x)) \, \mu(\mathrm{d}x) \right\} \\
  &= \min_{c \in \mathcal{C}^{(n)}(\bm{z})} \; \min_{r \in \mathcal{R}_c^{(n)}(\bm{z})} \; \int_\Xcal r(\bm{z} \bm{q}(x)) q(x) \, \mu(\mathrm{d}x) \notag \\
  &= \min_{c \in \mathcal{C}^{(n)}(\bm{z})} \; \min_{r \in \mathcal{R}_c^{(n)}(\bm{z})} \; \int_\Xcal r(\bm{z} \bm{q}(x)) \, Q_0(\mathrm{d}x) \\
  &= \min_{c \in \mathcal{C}^{(n)}(\bm{z})} \int_\Xcal c \, Q_0(\mathrm{d}x) \\
  &= \min \; \mathcal{C}^{(n)}(\bm{z}).
  \label{eq:d_derivative} 
\end{align}
Next, it is shown that the minimum of $\mathcal{C}^{(n)}(\bm{z})$ in \eqref{eq:d_derivative} can only take on values in $\partial_{z_0}^+ \rho_\lambda^{(n)}$. In order to see this, note that for any given $Q$ the set $\mathcal{C}^{(n)}(\bm{z})$ can be written as
\begin{equation}
  \mathcal{C}^{(n)}(\bm{z}) = \bigcap_{x \in \mathcal{S}_0} \partial_{z_0} \rho_\lambda^{(n)}(\bm{z} \bm{q}(x)),
\end{equation}
where $\mathcal{S}_0$ is the support of $Q_0$. Hence, $\mathcal{C}^{(n)}(\bm{z})$ is an intersection of superdifferentials, that is, of left-closed intervals on $\mathbb{R}_+$. If this intersection is empty, \eqref{eq:opt_condition_inequality} and \eqref{eq:opt_condition_equality} do not admit a solution, which contradicts the assumption that $Q_0$ is a maximizer. Hence, $\mathcal{C}^{(n)}(\bm{z})$ is the minimum of a non-empty intersection of left-closed intervals. It is not hard to show that this minimum is the largest left end-point of the intersecting intervals. Hence, there exits some $\tilde{x}$ such that
\begin{align}
  \partial_{z_0}^+ d_\lambda^{(n)}(\bm{z}) = \min \; \mathcal{C}^{(n)}(\bm{z}) &= \min \; \partial_{z_0} \rho_\lambda^{(n)}(\bm{z} \bm{q}(\tilde{x})) \\
  &= \partial_{z_0}^+ \rho_\lambda^{(n)}(\bm{z} \bm{q}(\tilde{x})) \\
  &\in \{0, \ldots, n\},
\end{align}
where the last step follows from the induction assumption. Since the arguments are independent of $\bm{z}$, it holds that
\begin{equation}
  \partial_{z_0}^+ d_\lambda^{(n)}(\bm{z}) \in \{0, \ldots, n\}
\end{equation}
for all $\bm{z} \in \mathbb{R}_+^3$. Using \eqref{eq:derivative_rho_n}, it follows that 
\begin{equation}
  \partial_{z_0}^+ \rho_\lambda^{(n+1)}(\bm{z}) \in \{0, \ldots, n+1\}.
\end{equation}
For $n \to \infty$, this yields
\begin{equation}
  \lim_{n \to \infty} \partial_{z_0}^+ \rho_\lambda^{(n+1)}(\bm{z}) = \partial_{z_0}^+ \rho_\lambda(\bm{z}) \in \Nbb_0 \cup \{\infty\}
\end{equation}
for all $\bm{z} \in \mathbb{R}_+^3$. Finally, the induction basis is given by $\partial_{z_0}^+ \rho_\lambda^{(0)} = \partial_{z_0}^+ g_\lambda = 0.$ 
This concludes the proof of the first part of the second statement.

In order to show the second part, consider the sets
\begin{align}
  \mathcal{Z}_0 = \left\{ (z_1, z_2) \in \mathbb{R}_+^2 : d_\lambda(0, z_1, z_2) = g_\lambda(z_1, z_2) \right\}, 
  \label{eq:z_zero} \\
  \mathcal{Z}_\infty = \left\{ (z_1, z_2) \in \mathbb{R}_+^2 : d_\lambda(0, z_1, z_2) < g_\lambda(z_1, z_2) \right\},
  \label{eq:z_infty}
\end{align}
with $d_\lambda$ defined in \eqref{eq:def_d}. Note that since $d_\lambda \leq \rho_\lambda$, $\{\mathcal{Z}_0, \mathcal{Z}_\infty \}$ is a partition of $\mathbb{R}_+^2$. By definition of $\Zcal_0$, it holds that
\begin{align}
  \rho_\lambda(\bm{z}) &= \min\{ g_\lambda(\bm{z}) \,,\, z_0 + d_\lambda(\bm{z}) \} = g_\lambda(\bm{z}).
  \label{eq:rho_g_equality}
\end{align}
for all $(z_1, z_2) \in \Zcal_0$ and all $z_0 \geq 0$. Therefore, it holds that
\begin{equation}
  \partial_{z_0}^+ \rho_\lambda(\bm{z}) = \partial_{z_0}^+ g_\lambda(\bm{z}) = 0
  \label{eq:drho_0}
\end{equation}
and, consequently, $\partial_{z_0}^+ \rho_\lambda(0, z_1, z_2) = 0$ for all $(z_1, z_2) \in \mathcal{Z}_0$. 

Now assume that $(z_1, z_2) \in \mathcal{Z}_\infty$. Since $d_\lambda = \lim_{n \to \infty} d_\lambda^{(n)}$, with $d_\lambda^{(n)}$ defined in \eqref{eq:def_dn}, this implies that there exists some $n^* \geq 1$ such that
\begin{equation}
  d_\lambda^{(n^*)}(0, z_1, z_2) < g_\lambda(z_1, z_2).
\end{equation}
Since $\bigl( d_\lambda^{(n)} \bigr)_{n \geq 0}$ is a nonincreasing sequence, this implies that 
\begin{align}
  \rho_\lambda^{(n)}(0, z_1, z_2) &= \min\left\{ g_\lambda(z_1, z_2) \,,\, d_\lambda^{(n)}(0, z_1, z_2) \right\} < g_\lambda(z_1, z_2)
\end{align}
for all $(z_1, z_2) \in \mathcal{Z}_\infty$ and all $n \geq n^*$. Hence, it follows from \eqref{eq:derivative_rho_n} that
\begin{equation}
  \partial_{z_0}^+ \rho_\lambda^{(n)}(0, z_1, z_2) = 1 + \partial_{z_0}^+ d_\lambda^{(n)}(0, z_1, z_2).
  \label{eq:drho_inf}
\end{equation}
for all $(z_1, z_2) \in \mathcal{Z}_\infty$ and all $n \geq n^*$. Since $d_\lambda^{(n)}$ is nondecreasing in $\bm{z}$ for all $n \geq 0$, it holds that $\partial_{z_0}^+ d_\lambda^{(n^*)} \geq 0$, so that 
\begin{equation}
  \partial_{z_0}^+ \rho_\lambda^{(n)}(0, z_1, z_2) \geq 1+n-n^*
  \label{eq:drho_bound}
\end{equation}
for all $(z_1, z_2) \in \mathcal{Z}_\infty$. Clearly, \eqref{eq:drho_bound} implies 
\begin{equation}
  \lim_{n \to \infty} \partial_{z_0}^+ \rho_\lambda^{(n)}(0, z_1, z_2) = \partial_{z_0}^+ \rho_\lambda(0, z_1, z_2) = \infty.
\end{equation}
This proves \eqref{eq:zero_sample_cost_drho}.

In order to prove \eqref{eq:zero_sample_cost_rho}, note that by definition of $\rho_\lambda$ in \eqref{eq:rho_implicit} it holds that that
\begin{align}
  \rho_\lambda(0, z_1, z_2) &= \min \left\{\, g_\lambda(\bm{z}) \,,\, \int_\Xcal \rho_\lambda\bigl(0, z_1 p_1(x), z_2 p_2(x)\bigr) \, \mu(\mathrm{d}x) \,\right\} \\
  &= \int_\Xcal \rho_\lambda\bigl(0, z_1 p_1(x), z_2 p_2(x)\bigr) \, \mu(\mathrm{d}x)
  \label{eq:rho0_implicit}
\end{align}
for all $(z_1, z_2) \in \mathbb{R}^2$, where the minimum is non-binding since $\rho_\lambda \leq g$. Since $\rho_\lambda$ is positively homogeneous and concave, it follows from a generalized version of Jensen's inequality \cite{Roselli2002} that in order for \eqref{eq:rho0_implicit} to hold, $\rho_\lambda(0, z_1, z_2)$ needs to be an affine function of $z_1, z_2$, that is,
\begin{equation}
  \rho_\lambda(0, z_1, z_2) = a z_1 + b z_2 + c
\end{equation}
for some $a, b, c \in \mathbb{R}_+$. However, since $\rho_\lambda$ is upper bounded by $g_\lambda$ it also needs to hold that
\begin{equation}
  a z_1 + b z_2 + c \leq g_\lambda(z_1, z_2) = \min\{ \lambda_1 z_1 \,,\, \lambda_2 z_2 \}
  \label{eq:rho_tilde_bound}
\end{equation}
for all $(z_1, z_2) \in \mathbb{R}_+^2$. It is not hard to show (consider the cases $z_1 = 0$ and $z_2 = 0$) that the only choice of $a, b, c$ that satisfies \eqref{eq:rho_tilde_bound} is $a = b = c = 0$, which in turn implies $\rho_\lambda(0, z_1, z_2) = 0$. This completes the proof.

\section{Proof of Theorem~\ref{th:optimality_sufficient}} 
\label{app:proof_th_3}                                    

Theorem~\ref{th:optimality_sufficient} defines the NP-KWT in terms of an optimal policy $\pi^*$ and a least favorable distribution $\Qbb_0$. A formal proof for the minimax optimality of this type of test is given in \cite{Fauss2020_aos} and will not be repeated here. However, given the specific choice of uncertainty sets in \eqref{eq:uncertainty_sets}, the optimal test can be stated in a more explicit manner.

The optimal decision and stopping rules in \eqref{eq:delta_opt} and \eqref{eq:psi_opt} are standard in sequential detection. The decision rule \eqref{eq:delta_opt} corresponds to a likelihood ratio test, and the stopping rule \eqref{eq:psi_opt} is obtained by comparing the cost of stopping with the expected cost of continuing; compare \cite{Novikov2009, Fauss2015}.

Sufficiency of the optimality conditions in \eqref{eq:lfd_product}, \eqref{eq:q_inequality} and \eqref{eq:q_equality} can be shown as follows. In \cite{Fauss2020_aos}, it is shown that in order for $\Qbb_0$ to be least favorable, it suffices that $\Qbb_0$ is of the form \eqref{eq:lfd_product}, and that
\begin{equation}
  Q_{\bm{z}} \in \argmax_{P_0 \in \Mcal_\mu} \; \int_{\Xcal} \rho_\lambda \bigl( \bm{z} \bm{p}(x) \bigr) \, \mu(\mathrm{d} x)
  \label{eq:Q_argmax}
\end{equation}
for all $\bm{z} \in \Rbb_+^3$. Using the same arguments as in the proof of Theorem~\ref{th:rho_properties}, a necessary and sufficient condition for $Q_{\bm{z}}$ to satisfy \eqref{eq:Q_argmax} is that there exists some $r \in \partial_{z_0} \rho_\lambda$ such that
\begin{align}
  c(\bm{z}) &\geq r \bigl( z_0 q_{\bm{z}}(x), z_1 p_1(x), z_2 p_2(x) \bigr) \quad \forall x \in \Xcal 
  \label{eq:q_argmax_inequality} \\
  c(\bm{z}) &= r \bigl( z_0 q_{\bm{z}}(x), z_1 p_1(x), z_2 p_2(x) \bigr) \quad Q_{\bm{z}}\text{-a.e.} \label{eq:q_argmax_equality}
\end{align}
for all $\bm{z} \in \Rbb_+^3$, where $c$ on the left-hand side is independent of $x$, but can depend on $\bm{z}$. Moreover, it follows from Theorem~\ref{th:rho_properties} that if \eqref{eq:q_argmax_inequality} and \eqref{eq:q_argmax_equality} admit a solution, $c(\bm{z})$ and $r$ can always be chosen such that $c(\bm{z}) \in \Nbb_0$. Note that $c(\bm{z}) = \infty$ is infeasible since it leads to the contradiction that $q_{\bm{z}} = 0$ needs to hold $Q_{\bm{z}}$-almost everywhere.

In order to connect the characterization of the least favorable distribution in \eqref{eq:q_argmax_inequality} and \eqref{eq:q_argmax_equality} to the optimal stopping policy, another result from \cite{Fauss2020_aos} is invoked. Namely, it is shown in \cite[Theorem~3]{Fauss2020_aos} that for policies $\pi^*$ with decision and stopping rules of the form \eqref{eq:delta_opt} and \eqref{eq:psi_opt}, respectively, it holds that
\begin{equation}
  \gamma_{\pi^*, \Qbb_0} \in \partial_{z_0} \rho_\lambda.
  \label{eq:sample_size_derivative}
\end{equation}
Substituting $\gamma_{\pi^*, \Qbb_0}$ for $r$ in \eqref{eq:q_argmax_inequality} and \eqref{eq:q_argmax_equality} yields \eqref{eq:q_inequality} and \eqref{eq:q_equality}. This completes the proof.


\section{Proof of Corollary~\ref{cl:threshold_test}} 
\label{app:proof_cl_1}                               

The statement in Corollary~\ref{cl:threshold_test} was implicitly already shown in the proof of Theorem~\ref{th:rho_properties}. The optimal stopping rule is obtained by comparing 
\begin{equation}
  z_0 + d_\lambda(\bm{z}) \lesseqqgtr g_\lambda(\bm{z}),
  \label{eq:certain_stopping}
\end{equation}
where the left-hand side corresponds to the expected cost for continuing the optimal test and the right hand side corresponds to the cost for stopping. Since $g_\lambda$ is independent of $z_0$ and the left-hand side is continuous, increasing, and unbounded in $z_0$, \eqref{eq:certain_stopping} can be rewritten as 
\begin{equation}
  z_0 \lesseqqgtr z_0^*(z_1, z_2),
\end{equation}
where $z_0^*(z_1, z_2)$ is the unique solution of
\begin{equation}
  z_0 + d_\lambda(\bm{z}) = g_\lambda(\bm{z}).
\end{equation}
This concludes the proof.

\section{Proof of Corollary~\ref{cl:lfd_support}} 
\label{app:proof_cl_2}                            

Let 
\begin{equation}
 \bm{z}\bm{q}(x) = \bigl( z_0 q(x),\, z_1 p_1(x),\, z_2 p_2(x) \bigr)
\end{equation}
with $Q \in \mathcal{Q}_{\bm{z}}$ and define the sets
\begin{align}
  \Xcal_0(\bm{z}) &\coloneqq \{ x \in \Xcal : (z_1 p_1(x), z_2 p_2(x)) \in \mathcal{Z}_0 \}, \\
  \Xcal_\infty(\bm{z}) &\coloneqq \{ x \in \Xcal : (z_1 p_1(x), z_2 p_2(x)) \in \mathcal{Z}_\infty \}.
\end{align} 
with $\mathcal{Z}_0$ and $\mathcal{Z}_\infty$ defined in \eqref{eq:z_zero} and \eqref{eq:z_infty}, respectively. 

In order to simplify the proof of the corollary, four auxiliary lemmas are stated first. All lemmas are simple consequences of already established results and are introduced mainly for the sake of a more orderly presentation. They are proven at the end of this section.
\begin{enumerate}
  \item $\rho_\lambda(\bm{z}) = 0$ for all $\bm{z} \ngtr 0$.
  \item $\Xcal_0(\bm{z}) = \bigl\{ x \in \Xcal : \min\{ p_1(x), p_2(x) \} = 0 \bigr\}$ for all $\bm{z} > 0$.
  \item $q(x) > 0$ for all $x \in \Xcal_\infty$.
  \item $Q(\Xcal_0(\bm{z})) > 0$ implies $\gamma_{\pi^*, \Qbb_0} \bigl( \bm{z}\bm{q}(x) \bigr) = 0$ for all $x \in \Xcal$
\end{enumerate}

The corollary can now be shown by considering two cases: For $\bm{z} \ngtr 0$, it follows from Lemma~1 that $\rho_\lambda(\bm{z}\bm{q}(x)) = 0$ for all $Q \in \Mcal_\mu$. Hence, $Q$ is arbitrary and the first statement in Corollary~\ref{cl:lfd_support} holds. For $\bm{z} > 0$, it follows from Lemma~2 and Lemma~3 that 
\begin{equation}
  \min\{ p_1(x), p_2(x) \} > 0 \quad \Rightarrow \quad q(x) > 0,
  \label{eq:Q_support}
\end{equation}
which implies the second statement in Corollary~\ref{cl:lfd_support}. Lemma~4 implies that $Q(\Xcal_0(\bm{z})) = 0$ whenever the test has a positive expected remaining sample size after having observed $X$. However, $Q\bigl(\min\{ p_1(X), p_2(X) \} = 0 \bigr) = 0$, in combination with \eqref{eq:Q_support}, implies $\supp(Q) = \supp(P_1) \cap \supp(P_2)$. Hence, the inclusion in the second statement can only be strict if the test stops after having observed $X$.

It remains to show the Lemmas.

\paragraph*{Lemma~1:} For $z_1 = 0$ and $z_2 = 0$, the lemma follows from $g_\lambda(\bm{z}) = \min\{ \lambda_1 z_1, \lambda_2 z_2\} = 0 \geq \rho_\lambda(\bm{z})$. For $z_0 = 0$, it follows from \eqref{eq:zero_sample_cost_rho} in Theorem~\ref{th:rho_properties}. \qed

\paragraph*{Lemma~2:} According to \eqref{eq:rho_g_equality}, it holds on $\Zcal_0$ that $\rho_\lambda(\bm{z}) = g_\lambda(\bm{z})$ for all $z_0 \geq 0$. However, according to \eqref{eq:zero_sample_cost_rho}, $\rho_\lambda(0, z_1, z_2) = 0$. Since $g_\lambda$ is independent of $z_0$, both statements can only be true if $\rho_\lambda(\bm{z}) = g_\lambda(\bm{z}) = \min\{ \lambda_1 z_1, \lambda_2 z_2 \} = 0$, which in turn implies $\min\{ z_1, z_2 \} = 0$. Consequently, $x \in \Xcal_0(\bm{z})$ implies $\min\{ z_1 p_1(x), z_2 p_2(x) \} = 0$. Since $\bm{z} > 0$ by assumption, this implies $\min\{ p_1(x), p_2(x) \} = 0$. \qed

\paragraph*{Lemma~3:} By definition, it holds that $\partial_{z_0}^+ \rho_\lambda\bigl(0, z_1, z_2) = \infty$ on $\Zcal_\infty$ . Since the partial differential at $z_0 = 0$ is unique, it follows from \eqref{eq:sample_size_derivative} that
\begin{equation}
  \gamma_{\pi^*, \Qbb_0} \bigl( \bm{z}\bm{q}(x) \bigr) = \infty,
  \label{eq:gamma_inf}
\end{equation}
for all $x \in \Xcal_\infty$ such that $q(x) = 0$. Clearly, \eqref{eq:gamma_inf} contradicts the optimality condition in \eqref{eq:q_inequality}, which states that the expected remaining sample sizes of all sequences need to be bounded away from infinity. \qed

\paragraph*{Lemma~4:} According to \eqref{eq:rho_g_equality}, it holds on $\Zcal_0$ that $\rho_\lambda(\bm{z}) = g_\lambda(\bm{z})$ for all $z_0 \geq 0$, which implies $\partial_{z_0}^+ \rho_\lambda(\bm{z}) = \partial_{z_0}^+ g_\lambda(\bm{z}) = 0$. Since the
partial differential at $z_0 = 0$ is unique, it follows from \eqref{eq:sample_size_derivative} that
\begin{equation}
  \gamma_{\pi^*, \Qbb_0} \bigl( \bm{z}\bm{q}(x) \bigr) = 0.
  \label{eq:gamma_zero}
\end{equation}
for all $x \in \Xcal_0$. Hence, according to the optimality conditions in \eqref{eq:q_inequality} and \eqref{eq:q_equality}, $Q(\Xcal_0) > 0$ implies $\gamma_{\pi^*, \Qbb_0} \bigl( \bm{z}\bm{q}(x) \bigr) = 0$ for all $x \in \Xcal$.

\section{Proof of Corollary~\ref{cl:truncation}} 
\label{app:proof_cl_3}                           

The sufficient conditon for an NP-KWT to be truncated is shown first. Let $\supp(P_1) \cap \supp(P_2) = \{ x^* \}$. By Corollary~\ref{cl:lfd_support}, this implies that $Q$ reduces to a single point mass at $x = x^*$ so that $z_0^n = 1$, $z_1^n = p_1^n(x^*)$, and $z_2^n = p_2^n(x^*)$ for all $n \geq 0$. That is, the NP-KWT becomes a regular KWT with $P_0 = \bm{1}_{x^*}$. It is not hard to show that this test is truncated at the smallest integer $n$ satisfing
\begin{equation}
  g_\lambda(z_1^{n-1}, z_2^{n-1}) - g_\lambda(z_1^n, z_2^n) \leq 1.
  \label{eq:kwt_truncation}
\end{equation} 
Note that for $\lambda_1 = \lambda_2 = \lambda$, \eqref{eq:kwt_truncation} can be solved explicitly, giving
\begin{equation}
  N = \left\lfloor 1 - \frac{\log \lambda + \log(1-p^*)}{\log p^*} \right\rfloor,
\end{equation}
where $p^* = \min\{p_1(x^*), p_2(x^*)\}$. 

For the proof of the sufficient conditions for an NP-KWT to be nontruncated, assume that 
\begin{equation}
  \biggl( \max_{x \in \supp(P_1) \cap \supp(P_2)} \frac{p_1(x)}{p_2(x)} \biggr) \biggl( \min_{x \in \supp(P_1) \cap \supp(P_2)} \frac{p_1(x)}{p_2(x)} \biggr) \geq 1.
  \label{eq:lr_maximum}
\end{equation}
In case the maximum and/or the minimum do not exist, it is instead assumed that for every $\varepsilon > 0$ there exists some $\tilde{x}$ such that
\begin{equation}
  P_2\biggl[ \frac{p_1(X)}{p_2(X)} \frac{p_1(\tilde{x})}{p_2(\tilde{x})} \leq 1 \biggr] < \varepsilon.
  \label{eq:lr_maximum_limit}
\end{equation}
If \eqref{eq:lr_maximum} or \eqref{eq:lr_maximum_limit} are not satisfied, $P_1$ and $P_2$ can simply be swapped.

The proof is based on an inductive argument. Assume that after the $n$th sample has been observed the test is in a state $\bm{z}^n$ in which it continues with a nonzero probability, that is,
\begin{equation}
  g_\lambda(\bm{z}^n) \geq z_0 + d_\lambda(\bm{z}^n),
  \label{eq:not_stopped_condition}
\end{equation}
with $d_\lambda$ is defined in \eqref{eq:def_d}. Moreover, assume that $\bm{z}^n$ is such that
\begin{equation}
  \lambda z_1^n \geq \lambda z_2^n.
  \label{eq:decision_inequality} 
\end{equation}
In order for the test to stop with probability one after the $(n+1)$th sample has been observed, it needs to hold that 
\begin{equation}
  \rho_\lambda(z_0^n q(x), z_1^n p_1(x), z_2^n p_2(x)) = g_\lambda(z_1^n p_1(x), z_2^n p_2(x))
  \label{eq:stopping_condition}
\end{equation}
for all $x \in \Xcal$. This follows directly from the optimality conditions in Theorem~\ref{th:optimality_necessary}. Integrating \eqref{eq:stopping_condition} yields
\begin{align}
  d_\lambda(\bm{z}^n) &= \int_\Xcal \rho_\lambda(z_0^n q(x), z_1^n p_1(x), z_2^n p_2(x)) \, \mu(\mathrm{d}x) \\
  &= \int_\Xcal g_\lambda(z_1^n p_1(x), z_2^n p_2(x)) \, \mu(\mathrm{d}x) \\
  &= \int_\Xcal \min\{ \lambda_1 z_1^n p_1(x), \lambda_2 z_2^n p_2(x) \} \, \mu(\mathrm{d}x).
\end{align}
Since
\begin{equation}
    \int_\Xcal \min\{ \lambda_1 z_1^n p_1(x), \lambda_2 z_2^n p_2(x) \} \, \mu(\mathrm{d}x) \leq \min\{\lambda_1 z_2^n, \lambda_2 z_2^n\} = \lambda_2 z_2^n,
\end{equation}
where the last equality follows from \eqref{eq:decision_inequality}, there exists some $\varepsilon_n \geq 0$ such that
\begin{equation}
  \int_\Xcal \min\{ \lambda_1 z_1^n p_1(x), \lambda_2 z_2^n p_2(x) \} \, \mu(\mathrm{d}x) = \lambda_2 z_2^n - \varepsilon_n.
  \label{eq:def_epsilon}
\end{equation}
Now, in order for \eqref{eq:not_stopped_condition} to be satisfied, it needs to hold that
\begin{align}
  g_\lambda(\bm{z}^n) &\geq z_0^n + d_\lambda(\bm{z}^n) \\
  \lambda_2 z_2^n &\geq z_0^n + \lambda_2 z_2^n - \varepsilon_n \\
   \varepsilon_n &\geq z_0^n
\end{align}
In turn, if $\varepsilon_n < z_0^n$ holds, \eqref{eq:stopping_condition} cannot be satisfied and the test is guaranteed to continue with positive probability after having observed the $(n+1)$th sample. 

Qualitatively speaking, saying that $\varepsilon_n$ is small is to say that the probability of changing the currently preferred hypothesis after having observed $X_{n+1}$ is small. The idea of the remainder of the proof is to show that it is possible to generate a sequence $\bm{x}^*$ such that the probability of any sample changing the preferred hypothesis becomes arbitrarily small. In a nutshell, this can be achieved by simply generating ``sufficiently significant'' observations.

In order to make this argument precise, define the sets
\begin{equation}
  \mathcal{C}(\varepsilon) \coloneqq \left\{ c \geq 0 : \int_\Xcal \min\{ c \, p_1(x), p_2(x) \} > 1-\varepsilon \right\}
\end{equation}
and
\begin{equation}
  \Xcal^*(t, \varepsilon) \coloneqq \left\{x \in \Xcal : \frac{p_1(x)}{p_2(x)} \geq 1, \; t \, \frac{p_1(x)}{p_2(x)} \in \mathcal{C}(\varepsilon) \right\}.
\end{equation}
By construction, it holds that any sequence $\bm{x}^*$ with 
\begin{equation}
  x_{n+1}^* \in \Xcal^*(t_n, \varepsilon_n), \quad n \geq 0,
  \label{eq:continuation_set}
\end{equation}
where
\begin{equation}
  t_n = \frac{\lambda_1 z_1^n}{\lambda_2 z_2^n} = t_{n-1}\frac{p_1(x_n^*)}{p_2(x_n^*)} \quad \text{and} \quad \varepsilon_n = \frac{z_0^n}{\lambda_2 z_2^n} = \varepsilon_{n-1} \frac{q(x_n^*)}{p_2(x_n^*)},
\end{equation}
is such that $\varepsilon_n < z_0^n$ for all $n \geq 1$. It remains to show that, under the assumptions in the corollary, the sets in \eqref{eq:continuation_set} are non-empty for all $n \geq 0$. First, in order for the set $\mathcal{C}(\varepsilon)$ to be nonempty for all $\varepsilon > 0$ it needs to hold that $\supp(P_1) = \supp(P_2)$, which is true by the assumption in the corollary. Second, given that $\mathcal{C}(\varepsilon)$ is nonempty, it follows from \eqref{eq:lr_maximum} or \eqref{eq:lr_maximum_limit} that the set $\Xcal^*(t, \varepsilon)$ is nonempty for all $t \geq 1$. Third, the assumption $\lambda_1 = \lambda_2$ guarantees that $t_0 = 1$ and, hence, by definition of $\Xcal^*(t, \varepsilon)$, that $t_n \geq 1$ for all $n \geq 0$. Finally, requiring the expected sample size to be larger than one guarantees that the test is not truncated after or even before having taken the first sample. This completes the proof.

\end{appendix}

\section*{Acknowledgements} 

The authors would like to sincerely thank Alexander Tartakovsky for the helpful conversations. His skepticism concerning our initial conjecture that the fixed-sample-size test solves the nonparametric Kiefer--Weiss problem motivated large parts of this work in the first place. 

This work was supported in part by the U.S.~National Science Foundation under Grant CCF-1908308.
The work of M.~Fau{\ss} was supported by the German Research Foundation (DFG) under Grant 424522268.

\bibliographystyle{plain}
\bibliography{nonpa_kiwei}

\end{document}